\documentclass[a4paper,11pt]{article}
\pdfoutput=1
\usepackage{srcltx}
\usepackage{indentfirst}
\usepackage{epsfig}
\usepackage{psfrag}
\usepackage{graphicx}
\usepackage{overpic}
\usepackage{multirow}
\usepackage{longtable}
\usepackage{hyperref}
\usepackage{subfig}
\usepackage{amsmath,amssymb, amsthm}
\usepackage{color}
\usepackage[mathscr]{euscript}
\usepackage{cases}
\usepackage{enumerate}
\usepackage[toc,page,title,titletoc,header]{appendix} 
\usepackage{cite}
\usepackage{rotating}

\usepackage{pgf,pgfarrows,pgfnodes,pgfautomata,pgfheaps}
\usepackage{tikz}
\usetikzlibrary{shapes.geometric} 
\usetikzlibrary{calc} 
\usetikzlibrary{backgrounds}
\usetikzlibrary{matrix} 
\usetikzlibrary{decorations.pathreplacing} 

\tikzset{
passprocess/.style={rectangle, minimum width=75pt, minimum height=25pt, draw=black},
startstop/.style={rectangle, rounded corners=5pt, minimum width=75pt, minimum height=25pt, draw=black},
decision/.style={
diamond,
shape aspect=3,
minimum width=75,
draw=black},
line/.style={draw, ->, 
thick
}} 

\newcommand{\bg}{\boldsymbol{g}}

\newcommand\bbR{\mathbb{R}}
\newcommand\bbN{\mathbb{N}}
\newcommand\bxi{\boldsymbol{\xi}}
\newcommand\bx{\boldsymbol{x}}

\newcommand\bq{\boldsymbol{q}}
\newcommand\bu{\boldsymbol{u}}
\newcommand\bv{\boldsymbol{v}}

\newcommand\bn{\boldsymbol{n}}

\newcommand\bLambda{{\boldsymbol{\Lambda}}}

\newcommand\bF{\boldsymbol{F}}

\newcommand\dd{\,\mathrm{d}}
\newcommand\He{\mathit{He}}
\newcommand\Kn{\mathit{Kn}}

\newcommand\mH{\mathcal{H}}

\newcommand\mM{\mathcal{M}}

\newcommand\NRxx{NR$xx$~}

\newcommand\mF{\mathcal{F}}

\newcommand\sss{\scriptscriptstyle}

\graphicspath{{images/}}
{\theoremstyle{remark} 
 \newtheorem{algorithm}{\bfseries Algorithm}
}

\title{An Efficient Steady-State Solver for Microflows with High-Order
  Moment Model}

\author{Zhicheng Hu\thanks{Department of Mathematics, College of
    Science, Nanjing University of Aeronautics and Astronautics,
    Nanjing 210016, China. Email: {\tt huzhicheng@nuaa.edu.cn}.}, ~~
  Guanghui Hu\thanks{Corresponding author. Department of Mathematics,
    University of Macau, Macao SAR, China; Zhuhai UM Science \& Technology
    Research Institute, Guangdong, China. Email: {\tt garyhu@umac.mo}.}}

\begin{document}

\maketitle

\begin{abstract}  
  In [Z. Hu, R. Li, and Z. Qiao. Acceleration for microflow
  simulations of high-order moment models by using lower-order model
  correction. J. Comput. Phys., 327:225-244, 2016], it has been
  successfully demonstrated that using lower-order moment model
  correction is a promising idea to accelerate the steady-state
  computation of high-order moment models of the Boltzmann
  equation. To develop the existing solver, the following aspects are
  studied in this paper. First, the finite volume method with linear
  reconstruction is employed for high-resolution spatial
  discretization so that the degrees of freedom in spatial space
  could be reduced remarkably without loss of accuracy. Second, by
  introducing an appropriate parameter $\tau$ in the correction step,
  it is found that the performance of the solver can be improved
  significantly, i.e., more levels would be involved in the solver,
  which further accelerates the convergence of the method. Third,
  Heun's method is employed as the smoother in each level to enhance
  the robustness of the solver. Numerical experiments in microflows
  are carried out to demonstrate the efficiency and to investigate the
  behavior of the new solver. In addition, several order reduction
  strategies for the choice of the order sequence of the solver are
  tested, and the strategy
  $m_{\sss l-1} = \lceil m_{\sss l} / 2 \rceil$ is found to be most
  efficient.

\vspace*{4mm}
\end{abstract}

\section{Introduction} 
\label{sec:intro} 

In the past few decades, the simulation of the Boltzmann equation has
attracted a great deal of attention in a variety of high-tech fields
such as rarefied gas dynamics in astronautics and fluid mechanics in
micro-electro-mechanical systems, where the mean free path of fluid
molecules becomes comparable to the characteristic length of the
problem. Because of the inherent high dimension of variables and the
complicated expression of the binary collision operator, an accurate
and efficient simulation of the Boltzmann equation still encounters
great challenges even for the computers nowadays. Lots of work has
been done to overcome these difficulties. One of the important efforts
is to reduce the computational cost of the collision operator by
employing simplified collision operators instead of the original one
\cite{BGK,Holway,Shakhov,harris2004introduction}, or developing fast
algorithms for it via spectral methods \cite{Hu2016,
  wang2018approximation}.

Another famous work is the Grad moment method first proposed in
\cite{Grad}, which tries to reduce the degrees of freedom in velocity
space without loss of accuracy by using a certain Hermite spectral
expansion with parameters adaptive to the local physical quantities of
the fluid. The derived system of equations is a semi-discretization of
the Boltzmann equation from numerical point of view, yet it is usually
regarded as the Grad moment model,
or macroscopic transport model in today's literature, see e.g.
\cite{Struchtrup}. This model is actually hierarchically extended with
respect to the expansion order, and is expected to converge to the
underlying Boltzmann equation rapidly as the expansion order
increases. Unfortunately, the original Grad moment models are found to
be lack of hyperbolicity \cite{Grad13toR13} and may yield unphysical
subshocks \cite{Grad1952}. A number of methods have been proposed to
regularize the Grad moment models \cite{McDonald, Levermore,
  Struchtrup2013, NRxx, Fan_new, framework, di2017linear}. Among them,
a systematic approach to guarantee the global hyperbolicity of the
moment model up to arbitrary order was introduced in
\cite{Fan_new,framework}, which makes the practical application of
high-order moment models possible. The resulting hyperbolic moment
models are of interest to us in the current paper.

In \cite{NRxx, Cai, Li, Microflows1D}, a systematic numerical method,
abbreviated as the \NRxx method, has been developed for the
regularized moment model of arbitrary order. The unified framework of
the \NRxx method makes the numerical implementation of the high-order
moment model without much difficulties. However, the developed
time-stepping \NRxx method turns out to be inefficient, when
steady-state simulations or models with a sufficiently large order are
taken into consideration. It can be seen in \cite{Microflows1D,
  wang2018approximation} that steady-state simulations of the moment
model with the order larger than $20$ may need to be applied for
numerical purpose. In such a situation, the moment model would include
thousands of nonlinear equations, which are deeply coupled with each
other. This immediately leads to an enormous amount of computational
cost, especially for the steady-state computation in which a long time
simulation is always required. Due to the importance of steady-state
simulations in microflows and the frequent employment of high-order
moment models, we are mainly concerned in this paper about the
acceleration of simulations in such cases.

Observing the fact that almost all equations of a moment model are
contained in the moment model with a larger order, it might be
possible to accelerate the computation of the high-order moment model
by using a lower-order moment model.
A natural way is to employ the solution of the lower-order moment
model to provide the initial guess in the computation of the
high-order moment model. Unfortunately, it is found from our
investigation that this approach does not help much in improving the
convergence of the simulation, although the convergence history would
become smoother. Inspired by the well-known multigrid method
\cite{brandt2011book, hackbusch1985book}, which could accelerate the
convergence of a basic iteration greatly by reducing error components
from the problem at various levels, it might be feasible to improve
the computational efficiency of high-order moment models by adopting a
lower-order moment model correction as the coarse grid
correction in multigrid method. Following the framework of nonlinear
multigrid method \cite{hackbusch1985book}, a nonlinear multi-level
moment (NMLM) solver for the high-order moment model could then be
obtained by providing appropriate transformation operators between the
moment models with different orders. Such an idea could be as
effective as expected also based on the following observation: the
resulting NMLM solver would not only be viewed as a multigrid solver
of velocity space for the Boltzmann equation, but also coincide to
some extent with the $p$-multigrid method
\cite{fidkowski2005pmg, helenbrook2008solving} or spectral multigrid
method \cite{ronquist1987spectral, maday1988spectral}, by recalling
the derivation of the moment model. In fact, this idea has been first
proposed and numerically investigated in our previous paper
\cite{hu2016acceleration}. To the best of our knowledge, it might be
the first effort on developing multigrid method of velocity space for
the Boltzmann equation. It is shown in \cite{hu2016acceleration} that
a significant improvement in efficiency of the steady-state
computation could be obtained even for the moment model with a
relatively small order such as $4$ and $5$, which indicates the idea
of using lower-order moment model correction is quite promising to
accelerate the simulation.

Although the solver in \cite{hu2016acceleration} worked well in the computation
of steady states of high-order moment models, there is still room left for
further improvement, from both the accuracy and the efficiency points of view.
First of all, since the piecewise constant approximation is used in the spatial
discretization, the numerical solution is of first order only, which is too
diffusive to deliver numerical solution with high resolution. Then from the
numerical experiments in \cite{hu2016acceleration}, it is found that the
stability of the solver is sensitive to the correction from the lower level,
i.e., the convergence of the solver will be negatively affected if the
correction from the lower level is directly used, while the situation can be
improved effectively by rescaling the correction. Furthermore, different
smoothing and order reduction strategies are tested in a variety of benchmark
problems, and numerical results highlight some insight on designing quality
method.

Based on the above consideration and observation, in this paper, we further
develop the solver proposed in \cite{hu2016acceleration}, from the following
aspects, 
\begin{itemize}
  \item The finite volume method with linear reconstruction is employed
    for spatial discretization of the target moment model, so that the
        degrees of freedom in spatial space could be reduced greatly while
        still being able to give accurate results in comparison to the
        first-order discretization which has been utilized in
        \cite{hu2016acceleration}. Following the basic idea of the \NRxx
        method, the derived discretization will have a unified form with
        respect to the order of the model, thus can also be solved under a
        unified framework for the moment model up to an arbitrary order.
      \item To enhance the stability of the resulting NMLM solver when a lot
        of levels are involved, a relaxation parameter is introduced in the
        step of updating the solution after each lower-order moment model
        correction is obtained. The computation of this correction step is
        also simplified a lot, so is much faster than the original way used in
        \cite{hu2016acceleration}.
      \item A second-order time-stepping scheme, namely, Heun's
        method, is used as the smoother of the NMLM solver. Based on our
        numeircal experience, there are several advantages by using Heun's
        method. Comparing to the SGS-Newton iteration proposed in
        \cite{hu2014nmg}, Heun's method can be implemented much easier, while
        comparing to the SGS-Richardson iteration proposed in
        \cite{hu2016acceleration}, Heun's method exhibits better performance,
        especially when a large Knudsen number is considered. It is worth
        mentioning that Heun's method would enhance the robustness of the
        NMLM solver.

      \item Numerical experiments of three benchmark spatially one-dimensional
        problems are carried out to investigate the performance and behavior
        of the new NMLM solver. Various order reduction strategies, including
        $m_{\sss l-1} = m_{\sss l}-1$, $m_{\sss l-1} = m_{\sss l} - 2$,
        $m_{\sss l-1} = m_{\sss l} - 4$, and
        $m_{\sss l-1} = \lceil m_{\sss l} / 2 \rceil$, are taken into account
        for the choice of the order sequence of the NMLM solver. It is shown
        that the convergence rate of the NMLM solver is effectively improved as the total
        levels increases. Among the order reduction strategies we have
        tested, it turns out that the best strategy is
        $m_{\sss l-1} = \lceil m_{\sss l} / 2 \rceil$.
\end{itemize}
The numerical experiments successfully show that both the numerical accuracy of
the solution and the computational efficiency of the solver are improved
significantly, compared with the ones in \cite{hu2016acceleration}.

The rest of this paper is arranged as follows. A brief review of the
underlying model equations in microflows as well as the related spatial
discretization is given in Section
\ref{sec:model}. The details of the nonlinear multi-level moment solver are
then described in Section \ref{sec:method}. Numerical experiments are
carried out in Section \ref{sec:example} to show the performance and
behavior of the proposed nonlinear multi-level moment solver. At last,
some concluding remarks are given in Section \ref{sec:conclusion}.


\section{The governing equations and their discretization}
\label{sec:model}
In this section, we briefly review the Boltzmann equation in steady
state, and the globally hyperbolic moment models,
then introduce a unified spatial discretization with linear
reconstruction for the given models.

\subsection{The steady-state Boltzmann equation}
\label{sec:model-boltzmann}
In the gas kinetic theory, the Boltzmann equation is used to describe
the evolution of gas molecules. It has the form
\begin{equation}
  \label{eq:boltzmann}
  \bxi \cdot \nabla_{\bx} f + \bF \cdot \nabla_{\bxi} f = Q(f),
\end{equation}
when the steady state of the fluid has been achieved. Here
$f(\bx, \bxi)$ is the molecular distribution function, in which
$\bx\in \Omega \subset \bbR^D$ $(D=1,2, \text{or } 3)$ and
$\bxi\in\bbR^{3}$ are the spatial position and the particle velocity
respectively. The vector $\bF$ stands for the acceleration of
molecules due to external force fields, and the right-hand side $Q(f)$
is the collision term. 
Upon the collision number assumption
(cf. \cite{Cowling,harris2004introduction}), it is given by
\begin{align}
  \label{eq:binary-collision}
  Q(f) = \int_{\bbR^{3}} \int_{\mathbb{S}_{+}^{2}}
  B(\vert \bxi-\bxi_*\vert, \bn ) (f'f_*'-ff_*) \dd \bn \dd \bxi_*,
\end{align}
where $f' = f(\bx, \bxi')$, $f_{*}=f(\bx,\bxi_*)$,
$f_{*}'= f(\bx, \bxi_*')$, and the pairs $(\bxi, \bxi_{*})$ and
$(\bxi', \bxi_*')$ are the pre- and post-collision velocities of a
colliding pairs of particles, with the unit vector
$\bn \in \mathbb{S}_{+}^{2}$ directed along the line joining the
centers of them. The collision kernel $B$ is a non-negative function
depending on the potential between gas molecules.

Such a binary collision term causes a great challenge in numerical
simulation. Simplified collision models, such as the BGK-type
relaxation models \cite{BGK,Holway,Shakhov} and the linearized
collision model \cite{harris2004introduction}, have been proposed to
replace it while still being able to predict the major physical
features of interest in a variety of situations.

The BGK-type collision term reads
\begin{align}
  \label{eq:col-relaxation}
  Q(f) = \nu (f^{\text{E}} - f), 
\end{align}
where $\nu$ is the average collision frequency assumed independent of
the particle velocity, and $f^{\text{E}}$ is the equilibrium
distribution function which depends on the specific choice of model:
\begin{itemize}
\item For the ES-BGK model \cite{Holway}, it is an anisotropic Gaussian
distribution defined by
\begin{equation}
  \label{eq:ES-f^E}
  f^{\text{E}}(\bx,\bxi) = \frac{\rho(\bx)}{m_* \sqrt{ \det[2\pi \bLambda(\bx)]}
  } \exp\left(-\frac{1}{2} ( \bxi - \bu(\bx) )^T [\bLambda(\bx)]^{-1}
    (\bxi-\bu(\bx)) \right),
\end{equation}
where $m_*$ is the mass of a single gas molecule, and
$\bLambda=(\lambda_{ij})$ is a $3\times 3$ matrix with
\begin{align*}
  \lambda_{ij}(\bx) = \theta(\bx) \delta_{ij} + \left( 1-
  \frac{1}{\Pr} \right) \frac{\sigma_{ij}(\bx)}{\rho(\bx)}, \quad i,
  j = 1,2,3, 
\end{align*}
in which $\delta_{ij}$ is the Kronecker delta symbol, and $\Pr$ is the
Prandtl number.
\item For the Shakhov model \cite{Shakhov}, it reads
\begin{align}
  \label{eq:Shakhov-f^E}
  f^{\text{E}}(\bx, \bxi) = \left[ 1+\frac{(1-\Pr)(\bxi-\bu(
    \bx))\cdot \bq( \bx) }{5\rho(\bx) [\theta(\bx)]^2} \left(
    \frac{ \vert \bxi-\bu(\bx) \vert^2}{\theta(\bx)} - 5
    \right)\right] f^{\text{M}}(
    \bx, \bxi),
\end{align}
where $f^M$ is the local Maxwellian given by 
\begin{align}
  \label{eq:BGK-f^E}
  f^{\text{M}}(\bx,\bxi) = \frac{\rho( \bx)}{m_* [2\pi
  \theta(\bx)]^{3/2}} \exp\left(-\frac{\vert \bxi -
  \bu(\bx) \vert^2}{2\theta(\bx)} \right).
\end{align}
\end{itemize}
In the above equations, $\rho$, $\bu$, $\theta$, $\sigma$, and $\bq$
are macroscopic physical quantities known as density, mean velocity,
temperature, stress tensor, and heat flux, respectively. They can be
computed from the distribution function $f$ as follows
\begin{align}
  \label{eq:moments}
\begin{aligned}
  & \rho(\bx) = m_* \int_{\bbR^3} f(\bx,\bxi) \dd \bxi, \quad
  \rho(\bx) \bu(\bx) = m_* \int_{\bbR^3} \bxi f(\bx,\bxi) \dd \bxi, \\
  & \rho(\bx) \vert \bu(\bx) \vert^2 + 3 \rho(\bx) \theta(\bx) = m_*
  \int_{\bbR^3} \vert \bxi \vert^2 f(\bx, \bxi) \dd \bxi, \\ &
  \sigma_{ij}(\bx) = m_* \int_{\bbR^3} (\xi_i - u_i(\bx))(\xi_j -
  u_j(\bx)) f(\bx,\bxi) \dd \bxi - \rho(\bx) \theta(\bx) \delta_{ij},
  \quad i,j = 1,2,3,\\ & \bq(\bx) = \frac{m_*}{2} \int_{\bbR^3} \vert
  \bxi-\bu(\bx) \vert^2 (\bxi-\bu(\bx)) f(\bx,\bxi) \dd \bxi.
\end{aligned}
\end{align}
It is noticed that when $\Pr =1$, both the ES-BGK model and the
Shakhov model reduce to the simplest BGK model \cite{BGK}, in which
$f^{\text{E}}$ is chosen as the local Maxwellian, i.e.,
$f^{\text{E}} \equiv f^{\text{M}}$.

In this paper, we adopt the BGK-type collision term as an example to
illustrate our algorithm. However, it is pointed out that the
framework of the present algorithm is also suitable for some other
collision models, as can be seen below.

\subsection{The moment model of order $M$}
\label{sec:model-HME}
To obtain the steady-state moment models for the Boltzmann equation
\eqref{eq:boltzmann}, we first expand the distribution function $f$
into a series
as
\begin{align}
  \label{eq:dis-expansion}
  f(\bx, \bxi) = \sum_{\alpha \in \bbN^{3}} f_\alpha(\bx)
  \mH_{\alpha}^{[\tilde{\bu}(\bx), \tilde{\theta}(\bx)]} (\bxi), 
\end{align}
where $f_{\alpha}(\bx)$ are the coefficients,
and 
$\mH_{\alpha}^{[\tilde{\bu}, \tilde{\theta}]}(\cdot)$
are the basis functions defined by
\begin{align}
  \label{eq:base}
  \mH_{\alpha}^{[\tilde{\bu}, \tilde{\theta}]}(\bxi) =
  \frac{1}{m_*(2\pi \tilde{\theta})^{^{3/2}}
  \tilde{\theta}^{^{|\alpha|/2}}} \prod\limits_{d=1}^3
  \He_{\alpha_d}({v}_d)\exp \left(-{{v}_d^2}/{2} \right),
  \quad {\bv} = \frac{\bxi-\tilde{\bu}}{
    \sqrt{\tilde{\theta}}},~\forall \bxi\in\bbR^3,
\end{align}
in which 
$|\alpha| = \alpha_1 + \alpha_2 + \alpha_3$, and $\He_n(\cdot)$ is the
Hermite polynomial of degree $n$, i.e.,
\begin{equation*}
  \He_n(x) = (-1)^n\exp \left( x^2/2 \right) \frac{\dd^n}{\dd x^n} 
  \exp \left(-x^2/2 \right).
\end{equation*}
The parameters $\tilde{\bu}$ and $\tilde{\theta}$ in the basis
functions are selected respectively as the local mean velocity $\bu$
and the local temperature $\theta$, which are determined from $f$
itself via \eqref{eq:moments}.
With this choice, we also have the following
relations
\begin{align}\label{eq:moments-relation}
  \begin{aligned}
    & f_0 = \rho, \qquad f_{e_1} = f_{e_2} = f_{e_3} = 0, \qquad
    \sum_{d=1}^3 f_{2e_d} = 0,
    \\ 
    & \sigma_{ij} = (1+\delta_{ij}) f_{e_i+e_j}, \quad q_i = 2
    f_{3e_i} + \sum_{d=1}^3 f_{2e_d+e_i}, \qquad i,j=1,2,3,
\end{aligned}
\end{align}
from \eqref{eq:moments}, where $e_1$, $e_2$, $e_3$ represent the
multi-indices $(1,0,0)$, $(0,1,0)$, $(0,0,1)$, respectively.

Based on the derivation of the globally hyperbolic moment system
proposed in \cite{Fan_new,Li,framework}, we then get a system of
equations for $\bu$, $\theta$, and $f_{\alpha}$, $|\alpha| \leq M$,
which is called the moment model of order $M$, as follows
\begin{equation}
  \label{eq:mnt-eqs}
  \begin{split}
    & \sum_{j=1}^D \Bigg[ \left( \theta \frac{\partial f_{\alpha -
          e_j}}{\partial x_j} + u_j \frac{\partial
        f_{\alpha}}{\partial x_j} + (1-\delta_{|\alpha|,M})(\alpha_j +
      1) \frac{\partial f_{\alpha+e_j}}{\partial x_j} \right) \\ &+
    \sum_{d=1}^3 \frac{\partial u_d}{\partial x_j} \left( \theta
      f_{\alpha-e_d-e_j} + u_j f_{\alpha-e_d} +
      (1-\delta_{|\alpha|,M}) (\alpha_j + 1) f_{\alpha-e_d+e_j}
    \right) \\ &+ \frac{1}{2} \frac{\partial \theta}{\partial x_j}
    \sum_{d=1}^3 \left( \theta f_{\alpha-2e_d-e_j} + u_j
      f_{\alpha-2e_d} + (1-\delta_{|\alpha|,M}) (\alpha_j + 1)
      f_{\alpha-2e_d+e_j} \right) \Bigg] \\ &= \sum_{d=1}^3 F_d
    f_{\alpha-e_d} + Q_{\alpha},
    \qquad |\alpha| \leq M,
    \end{split}
\end{equation}
where $F_d$ is the $d$th component of the acceleration $\bF$, and
$Q_{\alpha}$ are the coefficients in the expansion of the collision
term under the same basis functions as $f$, namely,
\begin{align}
  \label{eq:Q-expansion}
  Q(f) =  \sum_{\alpha \in \bbN^{3}}
  Q_{\alpha}(\bx) \mH_{\alpha}^{[\bu(\bx),\theta(\bx)]}  (\bxi).
\end{align}
For the BGK-type collision term \eqref{eq:col-relaxation}, we have
\begin{align*}
  Q_{\alpha} = \nu (f^{\text{E}}_\alpha - f_\alpha),
\end{align*}
where the analytical computational formula of $f^{\text{E}}_\alpha$
can be found in \cite{Li} and \cite{Microflows1D} for the Shakhov
model and the ES-BGK model respectively. For the binary collision
operator \eqref{eq:binary-collision} as well as the linearized
collision model \cite{harris2004introduction} with some special kernel
$B$, the computation of $Q_{\alpha}$ can be found in
\cite{wang2018approximation}.

Since the moment model \eqref{eq:mnt-eqs} contains the classic
hydrodynamic equations when $M\geq 2$, it is usually viewed as the
macroscopic transport model or the extended hydrodynamic model in the
literature. While from numerical point of view, it can be also viewed
as a semi-discretization of the Boltzmann equation in the velocity
space, by noting that the solution of it forms an approximation of the
distribution function by
\begin{align}
  \label{eq:truncated-dis}
  f(\bx, \bxi) \approx \sum_{\vert \alpha \vert \leq M} f_\alpha(\bx)
  \mH_{\alpha}^{[\bu(\bx), \theta(\bx)]} (\bxi).
\end{align}
This makes it much easier to develop numerical solvers for the moment
model of arbitrary order under a unified framework. Meanwhile, any
solver developed for the moment model can be also regarded as a solver
for the underlying Boltzmann equation.

Obviously, the moment model \eqref{eq:mnt-eqs} is a nonlinear system
coupling all moments, including the mean velocity $\bu$, the
temperature $\theta$, and the coefficients $f_{\alpha}$, together. And
it is easy to show that the number of equations in a moment model of
order $M$ is
\begin{align}
  \label{eq:total-moments-number}
\mM_M = {M+3 \choose 3}.
\end{align}
With the additional relations \eqref{eq:moments-relation}, we have
that the total number of independent variables is the same. It follows
that the system is very large, e.g., $\mM_{10} = 286$ and
$\mM_{26} = 3654$, resulting a huge computational cost for a general
designed numerical method, when a high-order moment model is taken
into account. However, a high-order moment model such as $M=10$ is
commonly employed in practical simulations, as can be seen in
\cite{Microflows1D, wang2018approximation}, where we can even see that
the moment model with $M=26$ or larger order is necessary for some
cases.

In the following, we use $g \in \mF_M^{[\tilde{\bu},\tilde{\theta}]}$ to denote
a truncated expression of a series similar to \eqref{eq:truncated-dis}, where
$\mF_M^{[\tilde{\bu},\tilde{\theta}]}$ is a linear space spanned by
$\mH_{\alpha}^{[\tilde{\bu}, \tilde{\theta}]}(\bxi)$ for all $\alpha$ with
$|\alpha| \leq M$.

\subsection{Spatial discretization with linear reconstruction}
\label{sec:model-dis}
From now on, we restrict ourselves to spatially one-dimensional case
for simplicity. Following the framework of the \NRxx method,
which was developed in \cite{NRxx, Cai, Li, Microflows1D}, we can
obtain a unified finite volume discretization for the moment model
\eqref{eq:mnt-eqs} of an arbitrary order. The main idea is to treat
all moments together as the truncated expansion
\eqref{eq:truncated-dis}, instead of dealing with them individually.

Suppose $\{x_{i}\}_{i=0}^N$ constitute a mesh of the spatial
domain $[0,L_{D}]$,
and $f_{i}(\bxi)$, $f_{i}^{L}(\bxi)$, and $f_{i}^{R}(\bxi)$ are the
discrete distribution function, respectively, on the center, the left
boundary, and the right boundary of the $i$th grid cell
$[x_{i},x_{i+1}]$. Then the finite volume discretization of the
Boltzmann equation \eqref{eq:boltzmann} over the $i$th cell reads
\begin{align}
  \label{eq:mnt-eq-dis}
  \frac{F(f_i^{R}(\bxi),f_{i+1}^{L}(\bxi)) - F(f_{i-1}^{R}(\bxi),
    f_i^{L}(\bxi))}{\Delta x_i} = G(f_i(\bxi)), 
\end{align}
where $\Delta x_i = x_{i+1} - x_{i}$ is the length of the $i$th cell,
$F(\cdot,\cdot)$ is the numerical flux defined at the boundaries of
the cell, and $G(\cdot)$ represents the discretization of the
acceleration and collision terms of the Boltzmann equation
\eqref{eq:boltzmann}. Let us further assume that
$f_{i}(\bxi) \in \mF_{M}^{[\bu_{i}, \theta_{i}]}$, that is,
\begin{align}
  \label{eq:dis-expansion-i}
  f_i(\bxi) = \sum_{|\alpha |\leq M} f_{i,\alpha}
  \mH_{\alpha}^{[{\bu}_i, {\theta}_i]}(\bxi),
\end{align}
where $\bu_{i}$ and $\theta_{i}$ are the local mean velocity and the
local temperature, respectively, such that the relation
\eqref{eq:moments-relation} holds for the coefficients
$f_{i,\alpha}$. Then by projecting all terms of \eqref{eq:mnt-eq-dis},
numerical fluxes $F(f_{i-1}^{R},f_i^{L})$, $F(f_i^{R},f_{i+1}^{L})$
and the right-hand side $G(f_i)$, into
$\mF_{M}^{[\bu_{i}, \theta_{i}]}$, and matching the resulting
coefficients in \eqref{eq:mnt-eq-dis} for the same basis function
$\mH_{\alpha}^{[{\bu}_i, {\theta}_i]}(\bxi)$, we can obtain a system
which equivalently is a discretization of the moment model
\eqref{eq:mnt-eqs} over the $i$th cell. Apparently, the set of
$\bu_{i}$, $\theta_{i}$ and $f_{i,\alpha}$ constitutes the solution of
the moment model \eqref{eq:mnt-eqs} on the $i$th cell. Consequently,
we would simply say $f_{i}(\bxi) \in \mF_{M}^{[\bu_{i}, \theta_{i}]}$
is the solution of the moment model on the $i$th cell below.

For the left boundary distribution function $f_{i}^{L}(\bxi)$ and the
right boundary distribution function $f_{i}^{R}(\bxi)$ of the $i$th
cell, which are assumed to belong to
$\mF_{M}^{[\bu_{i}^{L}, \theta_{i}^{L}]}$ and
$\mF_{M}^{[\bu_{i}^{R}, \theta_{i}^{R}]}$, respectively, it is enough
to give the computational formulae for parameters $\bu_{i}^{L}$,
$\bu_{i}^{R}$, $\theta_{i}^{L}$, $\theta_{i}^{R}$ and all expansion
coefficients $f_{i,\alpha}^{L}$, $f_{i,\alpha}^{R}$ with
$|\alpha|\leq M$. By linear reconstruction, they are calculated by
\begin{align}
  \label{eq:moments-recon}
  \begin{aligned}
  \bu_{i}^{L}= \bu_{i} - \frac{\Delta x_{i}}{2} \bg_{i}, & \qquad 
\bu_{i}^{R} = \bu_{i} + \frac{\Delta x_{i}}{2} \bg_{i}, \\ 
  \theta_{i}^{L}= \theta_{i} - \frac{\Delta x_{i}}{2} g_{i}, & \qquad
\theta_{i}^{R} = \theta_{i} + \frac{\Delta x_{i}}{2} g_{i}, \\
  f_{i,\alpha}^{L}= f_{i,\alpha} - \frac{\Delta x_{i}}{2} g_{i,\alpha}, & \qquad
f_{i,\alpha}^{R} = f_{i,\alpha} + \frac{\Delta x_{i}}{2} g_{i,\alpha}, 
\end{aligned}
\end{align}
where $\bg_{i}$, $g_{i}$ and $g_{i, \alpha}$ are reconstructed slopes
of the corresponding moments in the $i$th cell. A first-order
discretization can be obtained by setting all slopes to be $0$. While
in this paper we consider a second-order discretization by employing
\begin{align*}
  &\bg_{i} = \frac{\bu_{i+1} - \bu_{i-1}}{\Delta x_{i} + (\Delta
    x_{i-1} + \Delta x_{i+1}) / 2},  \\
  &g_{i} = \frac{\theta_{i+1} - \theta_{i-1}}{\Delta x_{i} + (\Delta
    x_{i-1} + \Delta x_{i+1}) / 2},  \\
  &g_{i,\alpha} = \frac{f_{i+1,\alpha} - f_{i-1,\alpha}}{\Delta x_{i} + (\Delta
    x_{i-1} + \Delta x_{i+1}) / 2}, 
\end{align*}
in which $\bu_{i\pm 1}$, $\theta_{i\pm 1}$ and $f_{i\pm 1, \alpha}$
are the solution of the moment model on the $(i\pm 1)$th cell.

Finally, from the explicit form of the moment model
\eqref{eq:mnt-eqs}, it is not difficult to deduce that the expansion
coefficients of $G(f_{i})$ in $\mF_{M}^{[\bu_{i}, \theta_{i}]}$ is
given by
$G_{i,\alpha} = \sum_{d=1}^3 F_{i,d} f_{i,\alpha-e_d} + Q_{i,\alpha}$.
Yet the calculation of the expansion coefficients of the numerical
flux $F(\cdot, \cdot)$ in $\mF_{M}^{[\bu_{i}, \theta_{i}]}$ is usually
required a transformation between two spaces, $\mF_M^{[\bu,\theta]}$
and $\mF_M^{[\tilde{\bu}, \tilde{\theta}]}$, since the function in
$\mF_{M}^{[\tilde{\bu}, \tilde{\theta}]}$ rather than
$\mF_{M}^{[\bu_{i}, \theta_{i}]}$ is always involved. Such a
transformation is the core of the \NRxx method, and has been provided
in \cite{NRxx, Qiao}. In our algorithm presented below, this
transformation will also be employed frequently without being
explicitly pointed out. Additionally, the numerical flux used in
\cite{Microflows1D} is adopted in our experiments for comparison.


\section{The nonlinear multi-level moment solver}
\label{sec:method}
This section is devoted to develop an efficient solver for a given
high-order moment model \eqref{eq:mnt-eqs} with the unified
second-order discretization \eqref{eq:mnt-eq-dis}, by using the
lower-order moment model correction. We first introduce a basic
iteration to solve the moment model of a certain order, then
illustrate the main ingredients of a nonlinear multi-level moment
solver for the high-order moment model.

\subsection{Basic iteration}
\label{sec:method-heun}
We would like to rewrite the discretization \eqref{eq:mnt-eq-dis} over
the $i$th cell into the form
\begin{align}
  \label{eq:residual-eq}
  R_i(f) = r_i(\bxi),
\end{align}
where $R_{i}(f)$ is the local residual on the $i$th cell given by
\begin{align}
  \label{eq:residual}
  R_i(f) = \frac{F(f_i^{R}(\bxi),f_{i+1}^{L}(\bxi)) -
  F(f_{i-1}^{R}(\bxi), f_i^{L}(\bxi))}{\Delta x_i} - G(f_i(\bxi)),
\end{align}
and $r_i(\bxi) \in \mF_M^{[\bu_i, \theta_i]}$ is a known function
introduced to make \eqref{eq:residual-eq} suitable for a slightly more
general problem. For the discretization \eqref{eq:mnt-eq-dis}, we have
$r_i(\bxi) \equiv 0$. It is clear that the above discretization gives
a nonlinear system coupling all unknowns, i.e., $\bu_i$, $\theta_i$
and $f_{i,\alpha}$, with $i=0,1,\ldots,N-1$, and $|\alpha | \leq M$,
together. As stated in \cite{hu2016acceleration}, it is quite
difficult to design an efficient iteration for such a nonlinear system
based on the Newton-type method, especially for the case that the
order $M$ is sufficiently large. Alternatively, a simple relaxation
method, referred to the SGS-Richardson iteration, was proposed in
\cite{hu2016acceleration,hu2015} for the discretization
\eqref{eq:residual-eq} without linear reconstruction. It turns out
that the SGS-Richardson iteration could also work for the second-order
discretization \eqref{eq:residual-eq}. Nevertheless, we would employ
Heun's method instead of it in the current implementation for better
performance in the situation when the acceleration by using
lower-order moment model correction is considered.

Given an approximate solution $f_i^n(\bxi)$, $i=0,1,\ldots, N-1$,
Heun's method
first calculates an intermediate approximation $f_i^*(\bxi)$,
$i=0,1,\ldots, N-1$, by
\begin{align}
  \label{eq:heun-step-1}
  f_i^{*}(\bxi) = f_i^{n}(\bxi) + \omega \left(r_i(\bxi) -
    R_i(f^{n}) \right),
\end{align}
and then get the new approximate solution $f_i^{n+1}(\bxi)$,
$i=0,1,\ldots, N-1$, by
\begin{align}
  \label{eq:heun-step-2}
  f_i^{n+1}(\bxi) = f_i^{n}(\bxi) + \omega \left( r_i(\bxi) -
    \frac{1}{2} \left(R_i(f^{n}) + R_i(f^{*})\right) \right),
\end{align}
where the parameter $\omega$ is selected according to the CFL
condition
\begin{align*}
  \omega \max_{i} \left\{ \frac{\lambda_{\max,i}}{\Delta x_i} \right\} < 1,
\end{align*}
in which $\lambda_{\max,i}$ is the largest value among the absolute
values of all eigenvalues of the hyperbolic moment model
\eqref{eq:mnt-eqs} on the $i$th cell. Similar to the SGS-Richardson
iteration, each calculation of \eqref{eq:heun-step-1} and
\eqref{eq:heun-step-2} does numerically consist of two steps. As an
example, for \eqref{eq:heun-step-1}, we first find an approximation
$f_i^{**}(\bxi)$ in $\mF_M^{[\bu_i^{n},\theta_i^{n}]}$, such that its
expansion coefficients $f_{i,\alpha}^{**}$ in terms of the basis
functions $\mH_{\alpha}^{[\bu_i^n, \theta_i^n]} (\bxi)$ are calculated
by
\begin{align*}
  f_{i,\alpha}^{**} = f_{i,\alpha}^{n} + \omega \left(r_{i,\alpha} -
  R_{i,\alpha}\right), \quad |\alpha | \leq M,
\end{align*}
where $f_{i,\alpha}^{n}$, $r_{i,\alpha}$, and $R_{i,\alpha}$ represent
expansion coefficients respectively of $f_{i}^{n}(\bxi)$, $r_i(\bxi)$
and $R_i(f^{n})$ in terms of the same basis functions. Then we calculate
$\bu_{i}^{*}$ and $\theta_{i}^{*}$ from $f_{i}^{**}(\bxi)$, and
project $f_i^{**}(\bxi)$ into $\mF_M^{[\bu_i^{*},\theta_i^{*}]}$ to
obtain $f_i^{*}(\bxi)$.

A single level solver, for the moment model \eqref{eq:mnt-eqs} of a
certain order on a given mesh, is then obtained by performing Heun's
method repeatedly until the norm of the global residual $\tilde{R}$
with $\tilde{R}_i(\bxi) = r_i(\bxi) - R_i(f)$ is smaller than a given
tolerance, which indicates the steady state has been achieved. Here,
the same norm as in \cite{hu2016acceleration} is adopted in our
numerical experiments.


\subsection{Lower-order moment model correction}
\label{sec:method-lower}
In order to improve the efficiency of steady-state computation when
the moment model with a high order $M$ is involved, we now turn to
consider the acceleration strategy using the lower-order moment model
correction, as proposed in \cite{hu2016acceleration}. The key point is
to establish an appropriate relationship between the high-order problem
and the lower-order problem.

For convenience, the underlying problem resulting from the
discretization \eqref{eq:residual-eq} of a high order $M$ is rewritten
into a global form as
\begin{align}
  \label{eq:high-order-problem}
  R_M(f_M) = r_{\sss M}.
\end{align}
Suppose we get an approximate solution for the above problem and
denote it by $\bar{f}_{M}$ with its $i$th component
$\bar{f}_{M,i}(\bxi) \in \mF_M^{[\bar{\bu}_{M, i},
  \bar{\theta}_{M,i}]}$. Then following \cite{hu2016acceleration}, the
lower-order problem can be defined by
\begin{align}
  \label{eq:lower-order-problem}
  R_m(f_m) = r_m \triangleq R_m(\tilde{I}_M^m \bar{f}_M) + I_M^m \left(r_{\sss M} -
    R_M(\bar{f}_M)\right),
\end{align}
where $I_{M}^{m}$ and $\tilde{I}_M^m$ are the restriction operators
transferring functions from the high $M$th-order function space into a
lower $m$th-order function space, and usually do not require the same.
The lower-order operator $R_m$ is the same discretization operator as
the high-order counterpart $R_M$, except that $R_m$ is applied on the
moment model of a lower order $m$. As a result, the lower-order
problem \eqref{eq:lower-order-problem} can be solved by the same
method as the high-order problem \eqref{eq:high-order-problem}. Once
the solution $f_m$ of the lower-order problem
\eqref{eq:lower-order-problem} is obtained, the solution of the
high-order problem \eqref{eq:high-order-problem} could be then
corrected by
\begin{align}
  \label{eq:update-high-order}
  \hat{f}_M = \bar{f}_M + \tau I_m^M\left( f_m - \tilde{I}_M^m \bar{f}_M
  \right),
\end{align}
where $I_m^M$ is the prolongation operator transferring functions from
the $m$th-order function space to the $M$th-order function space, and
$\tau \in (0,1]$ is a relaxation parameter introduced to enhance the
stability of the final solver. For the case $\tau=1$, it reduces to
the correction employed in \cite{hu2016acceleration}.

\subsection{Restriction and prolongation}
\label{sec:method-restriction}
Currently, only the case that both high-order problem
\eqref{eq:high-order-problem} and lower-order problem
\eqref{eq:lower-order-problem} are defined on the same spatial mesh is
taken into consideration. For such a case, it is sufficient to give the
definition of the restriction and prolongation operators on an
individual element of the spatial mesh. Hence, we omit the index $i$
of the spatial element below without causing confusion.

It is not easy to design proper restriction and prolongation operators
directly based on the high-order moment set
$\{\bu_{M}, \theta_{M}, f_{M,\alpha}, |\alpha | \leq M \}$, and the
lower-order one 
$\{\bu_{m}, \theta_{m}, f_{m,\alpha}, |\alpha | \leq m \}$. With the
help of the unified expression \eqref{eq:dis-expansion-i} combining
all moments together, however, these transferring operators could be
constructed and implemented very simple and efficient following the
idea of the $p$-multigrid method \cite{fidkowski2005pmg,
  helenbrook2008solving}.

Note that the high-order solution $\bar{f}_{M}$ as well as the
associated residual $r_{\sss M} - R_M(\bar{f}_M)$ are expressed as a
function in $\mF_M^{[\bar{\bu}_M, \bar{\theta}_M]}$. And the initial
discretization of the lower-order problem
\eqref{eq:lower-order-problem} is formulated in
$\mF_m^{[\bar{\bu}_{m}, \bar{\theta}_{m}]}$ with
$\bar{\bu}_m = \bar{\bu}_{\sss M}$ and
$\bar{\theta}_m = \bar{\theta}_M$, as explained in
\cite{hu2016acceleration}. It follows that the basis functions of
$\mF_m^{[\bar{\bu}_{m}, \bar{\theta}_{m}]}$ coincide with the first
$\mM_{m}$ functions of the basis functions of
$\mF_M^{[\bar{\bu}_M, \bar{\theta}_M]}$. Using the orthogonality of
the basis functions, we thus define both the solution restriction
operator $\tilde{I}_M^m$ and the residual restriction operator
$I_{M}^{m}$ as the truncation operator that simply gets rid of the
part in terms of the basis functions 
$\mH_{\alpha}^{[\bar{\bu}_M, \bar{\theta}_M]}(\bxi)$ with
$|\alpha| > m$.

In the following, we will give the implementation of the correction step used
in this paper, in which the prolongation operator will also be described in
detail. First of all, the correction step in \cite{hu2016acceleration} can be
summarized as
\begin{enumerate}
\item Compute the lower-order correction
  $f_m - \tilde{I}_M^m \bar{f}_M$ in
  $\mF_m^{[\bar{\bu}_M, \bar{\theta}_M]}$, by calling the
  transformation from $\mF_m^{[{\bu}_m, {\theta}_m]}$ into
  $\mF_m^{[\bar{\bu}_M, \bar{\theta}_M]}$. 
\item Retruncate the lower-order correction from
  $\mF_m^{[\bar{\bu}_M, \bar{\theta}_M]}$ into
  $\mF_M^{[\bar{\bu}_M, \bar{\theta}_M]}$, by remaining the
  coefficients 
  with $|\alpha | \leq m$ unchanged, and setting the coefficients
  with $|\alpha | > m$ to be $0$.
\item Add the right-hand side of \eqref{eq:update-high-order} together
  in $\mF_M^{[\bar{\bu}_M, \bar{\theta}_M]}$, compute the new
  macroscopic velocity $\hat{\bu}_{M}$ and temperature
  $\hat{\theta}_{M}$, and then get the new approximate solution
  $\hat{f}_{M}$ by calling the transformation from
  $\mF_M^{[\bar{\bu}_M, \bar{\theta}_M]}$ into
  $\mF_M^{[\hat{\bu}_M, \hat{\theta}_M]}$.
\end{enumerate}

Instead of the above three steps in \cite{hu2016acceleration}, we propose the
following strategy for the correction
\begin{align}
  \label{eq:updated-high-order-moments}
  \begin{aligned}
    & \hat{\bu}_{M} = (1-\tau) \bar{\bu}_{M} + \tau \bu_{m}, \quad \hat{\theta}_{M}
  = (1-\tau) \bar{\theta}_{M} + \tau \theta_{m}, \\
  & \hat{f}_{M,\alpha} = \left\{
  \begin{aligned}
    & (1-\tau) \bar{f}_{M,\alpha} + \tau f_{m, \alpha}, && \quad |\alpha | \leq m, \\
    & \bar{f}_{M, \alpha}, && \quad m < |\alpha |\leq M.
  \end{aligned}\right.
  \end{aligned}
\end{align}

The motivation for proposing the new correction is based on the following
observation.

For the case $\tau=1$ in \cite{hu2016acceleration}, it is found 
that the coefficients of the right-hand side
\eqref{eq:update-high-order}, corresponding to the basis functions
$\mH_{\alpha}^{[\bar{\bu}_M, \bar{\theta}_M]}(\bxi)$ with
$|\alpha | \leq m$, entirely come from the projection of $f_{m}$ in
$\mF_M^{[\bar{\bu}_M, \bar{\theta}_M]}$, since we have
\begin{align*}
  \bar{f}_{M} - I_{m}^{M} \tilde{I}_{M}^{m} \bar{f}_{M} =
  \sum_{m<|\alpha |\leq M} \bar{f}_{M,\alpha}
  \mH_{\alpha}^{[\bar{\bu}_M, \bar{\theta}_M]}(\bxi).
\end{align*}
As the macroscopic velocity $\hat{\bu}_{M}$ and temperature
$\hat{\theta}_{M}$ depend only on the $\bar{\bu}_{M}$, $\bar{\theta}_{M}$
and expansion coefficients with $|\alpha | \leq 2$ from
\eqref{eq:moments} and \eqref{eq:dis-expansion}, we can deduce that
$\hat{\bu}_{M} = \bu_{m}$ and $\hat{\theta}_{M} = \theta_{m}$ if the
lower order $m\geq 2$. Therefore, the correction step can be
implemented more efficient by avoiding the transformation between
different function spaces, and simply setting $\hat{f}_{M}$ as
\begin{align*}
  \hat{\bu}_{M} = \bu_{m}, \quad \hat{\theta}_{M} = \theta_{m}, \quad
  \hat{f}_{M,\alpha} = \left\{
  \begin{aligned}
    & f_{m, \alpha}, && |\alpha | \leq m, \\
    & \bar{f}_{M, \alpha}, && m < |\alpha |\leq M.
  \end{aligned}\right.
\end{align*}
Although the above $\hat{f}_{M, \alpha}$ with $m < |\alpha |\leq M$ is
slightly different from the previous calculation, the performance of
the final solver is similar in our numerical experiments. 

It is noted that only $\tau = 1$ can be handled in \cite{hu2016acceleration},
while our new correction strategy can be applied to the more general
cases when $\tau\neq 1$.


\subsection{Multi-level algorithm}
\label{sec:method-velocity}
If the lower-order problem \eqref{eq:lower-order-problem} still has a
relatively large order, it is straightforward to solve it by employing
a much lower-order moment model correction as illustrated in previous
subsections. A nonlinear multi-level moment (NMLM) iteration for the
underlying discretization problem \eqref{eq:mnt-eq-dis} is then
obtained by recursively applying the lower-order moment model
correction.

Let us introduce $m_{\sss l}$, $l=0,1,\ldots,L$, to denote the order
of the $l$th-level problem, and suppose
$2\leq m_0 < m_1 <\cdots < m_{\sss L}$. Then the $(l+1)$-level NMLM
iteration produces the new approximate solution $f_{m_l}^{n+1}$ from a
given approximate solution $f_{m_l}^{n}$, denoted by
$f_{m_l}^{n+1} = \text{NMLM}_l(f_{m_l}^n, r_{m_l})$, as the following
algorithm.
\begin{algorithm} [Nonlinear multi-level moment (NMLM) iteration]~
\label{alg:nmg}
  \begin{enumerate}
  \item If $l=0$, call the lowest-order solver, which would be given
    later, to obtain the new approximate solution $f_{m_0}^{n+1}$; 
    otherwise, go to the next step.
  \item Pre-smoothing: perform $s_1$ steps of Heun's method beginning
    with the initial approximation $f_{m_l}^n$ to obtain a new
    approximation $\bar{f}_{m_l}$.
  \item Lower-order moment model correction:
    \begin{enumerate}
    \item Compute the high-order residual as $\bar{R}_{m_l} = r_{m_l}
      - R_{m_l}(\bar{f}_{m_l})$.
    \item Prepare the initial guess of the lower-order problem by the
      restriction operator $\tilde{I}_{m_l}^{m_{l-1}}$ as
      $\bar{f}_{m_{l-1}} = \tilde{I}_{m_l}^{m_{l-1}} \bar{f}_{m_l}$.
    \item Calculate the right-hand side of the lower-order problem
      \eqref{eq:lower-order-problem} as $r_{m_{l-1}} =
      I_{m_l}^{m_{l-1}} \bar{R}_{m_l} + R_{m_{l-1}}(
      \bar{f}_{m_{l-1}})$.
    \item Recursively call the NMLM iteration (repeat $\gamma$ times
      with $\gamma=1$ for a $V$-cycle, $\gamma=2$ for a
      $W$-cycle, and so on) to get the new approximation of the
      lower-order problem as
      \begin{align*}
        \tilde{f}_{m_{l-1}} = \text{NMLM}_{l-1}^\gamma (\bar{f}_{m_{l-1}},
        r_{m_{l-1}}).
      \end{align*}
    \item Correct the high-order solution $\hat{f}_{m_{l}}$ by the
      formula \eqref{eq:updated-high-order-moments}.
    \end{enumerate}
  \item Post-smoothing: perform $s_2$ steps of Heun's method beginning
    with $\hat{f}_{m_{l}}$ to obtain the new approximation
    $f_{m_l}^{n+1}$.
  \end{enumerate}
\end{algorithm}
Performing the above $(L+1)$-level NMLM iteration until the steady
state has been achieved, we consequently get an $(L+1)$-level NMLM
solver for the moment model \eqref{eq:mnt-eqs} of order $m_{\sss L}$.
Obviously, the one-level NMLM solver reduces to the single level
solver of Heun's method.

For the lowest-order solver, the Heun method is applied again by
noting that the lowest-order problem is analogous to the problem on
the other order levels. As the lowest-order problem is indeed not
necessary to be solved accurately, only $s_3$ steps of Heun's method
will be performed in each calling of the lowest-order solver, to make
the final NMLM solver more efficient. Here $s_3$ is a positive integer
a little larger than the smoothing steps $s_1+s_2$.

Another technical issue is the choice of the order sequence
$m_{\sss l}$, $l=L,L-1,\ldots,1,0$, of the NMLM solver. Three order
reduction strategies, i.e., $m_{\sss l-1} = m_{\sss l}-1$,
$m_{\sss l-1} = m_{\sss l} - 2$, and
$m_{\sss l-1} = \lceil m_{\sss l} / 2 \rceil$, have been numerically
investigated in \cite{hu2016acceleration}. It turns out that all three
order reduction strategies could usually accelerate the steady-state
computation, and the most efficient strategy should be
$m_{\sss l-1} = \lceil m_{\sss l} / 2 \rceil$, the second should be
$m_{\sss l-1} = m_{\sss l} - 2$, and the third should be
$m_{\sss l-1} = m_{\sss l} - 1$. In the next section, we will
investigate the performance of all these three order reduction
strategies again on the proposed new NMLM solver.

From our observation, the convergence rate and the efficiency of the
NMLM solver usually become better as the total levels
increases. However, the choice $\tau=1$ in the correction step
\eqref{eq:updated-high-order-moments} to some extent introduces
instability of the NMLM solver when too many levels is employed, and a
too small $\tau$, e.g., $\tau = 0.5$ would also make the NMLM solver
inefficient. As the optimal $\tau$ is difficult to be determined and
need to be further studied, we currently set $\tau=0.9$ throughout our
numerical experiments.


\section{Numerical experiments}
\label{sec:example}

Three numerical examples, i.e., the planar Couette flow, the force
driven Poiseuille flow and the Fourier flow, are given in this section
to illustrate the main features of the proposed NMLM solver. The
dimensionless case with the molecular mass $m_*$ to be 1 is considered
without loss of generality. To complete the problem, the Maxwell
boundary conditions derived in \cite{Li} for the moment model are
employed for all examples. Since such boundary conditions could not
determine a unique solution for the steady-state moment model
\eqref{eq:mnt-eqs} as mentioned in \cite{hu2014nmg}, the solution is
corrected as in \cite{Mieussens2004, hu2014nmg} at each NMLM
iteration, to recover the consistent steady-state solution with the
time-stepping scheme and the solvers proposed in \cite{hu2014nmg} and
\cite{hu2016acceleration}.

In all numerical tests, the $V$-cycle NMLM solver with $s_1=s_2=2$ and
$s_3=5$ is performed under CentOS system on an Xeon
workstation with a 12-core processor and core speed 3.00GHz. All
computations are starting at the global Maxwellian with
\begin{align}
  \label{eq:initial}
  \rho^0 (x) = 1, \quad \bu^0(x) = 0, \quad \theta^0(x) = 1.
\end{align}
The tolerance indicating the achievement of steady state is set to be
$10^{-8}$. We have observed that the behavior of the NMLM solver are
very similar for the BGK-type collision models. Hence only results for
the ES-BGK collision model with the Prandtl number $\Pr=2/3$ are
presented below.

\subsection{The planar Couette flow}
\label{sec:num-ex-couette}
We first consider the planar Couette flow, a frequently used benchmark
test in microflows. The same settings as in \cite{Microflows1D,
  hu2014nmg} are adopted in our tests. To be specific, the gas of argon
is considered in the space between two infinite parallel plates that
is separated by a distance of $L_{D}=1$. Both plates have the
temperature $\theta^W=1$, and move in the opposite direction along the
plate with a relative speed $u^{W} = 1.2577$. The dimensionless
collision frequency $\nu$ is given by
\begin{align}
  \label{eq:couette-nu}
  \nu = \sqrt{\frac{\pi}{2}} \frac{\Pr}{\Kn} \rho \theta^{1-w},
\end{align}
where $\Kn$ is the Knudsen number, and $w$ is the viscosity index set
to be $0.81$. There is no external force acting on the gas, i.e.,
$\bF\equiv 0$, so the gas is only driven by the motion of the plates
and would finally reach a steady state.

\begin{figure}[!htb]
  \centering
  \subfloat[Density, $\rho$]{\includegraphics[width=0.49\textwidth]{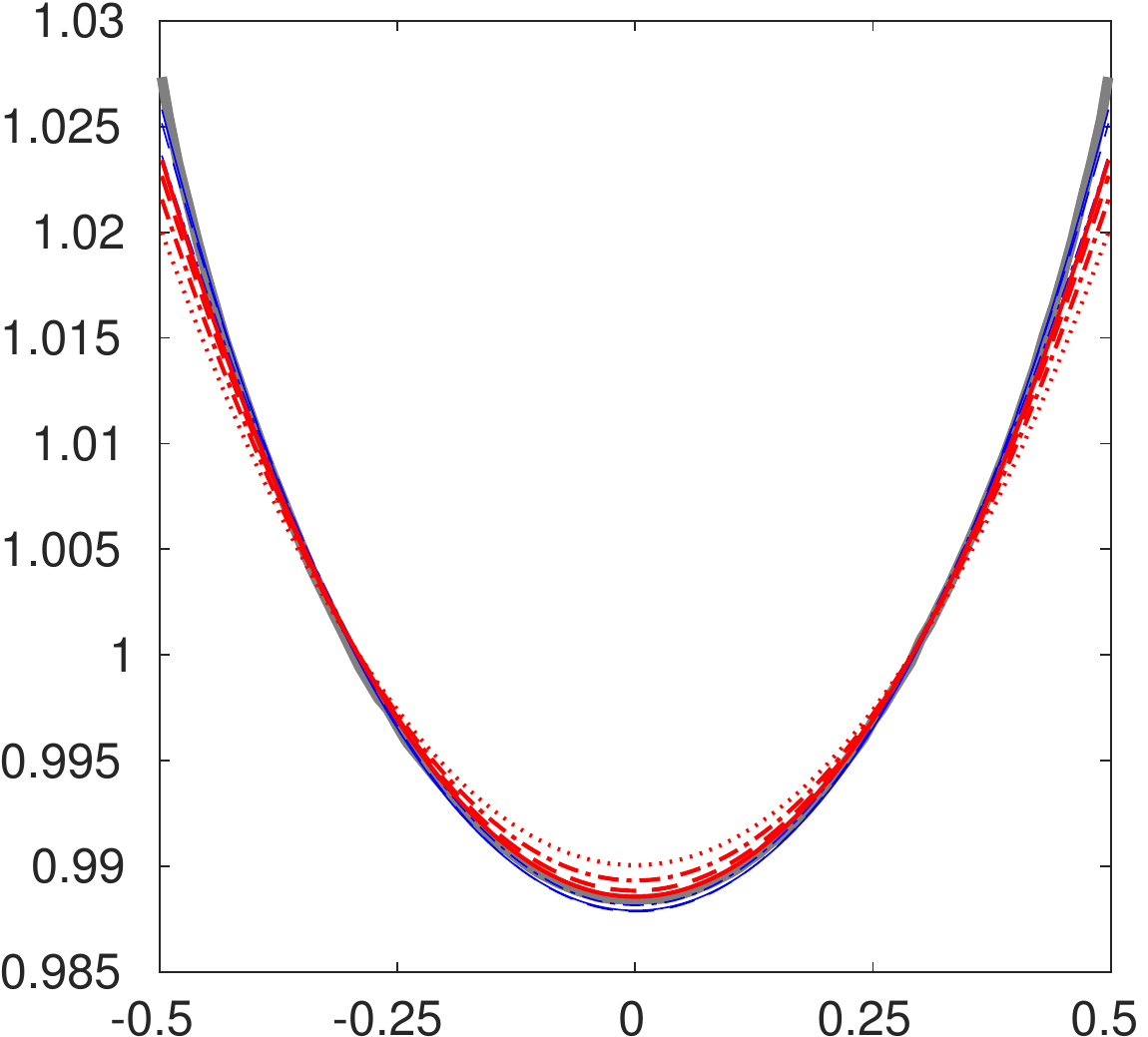}}\hfill
  \subfloat[Temperature, $\theta$]{\includegraphics[width=0.49\textwidth]{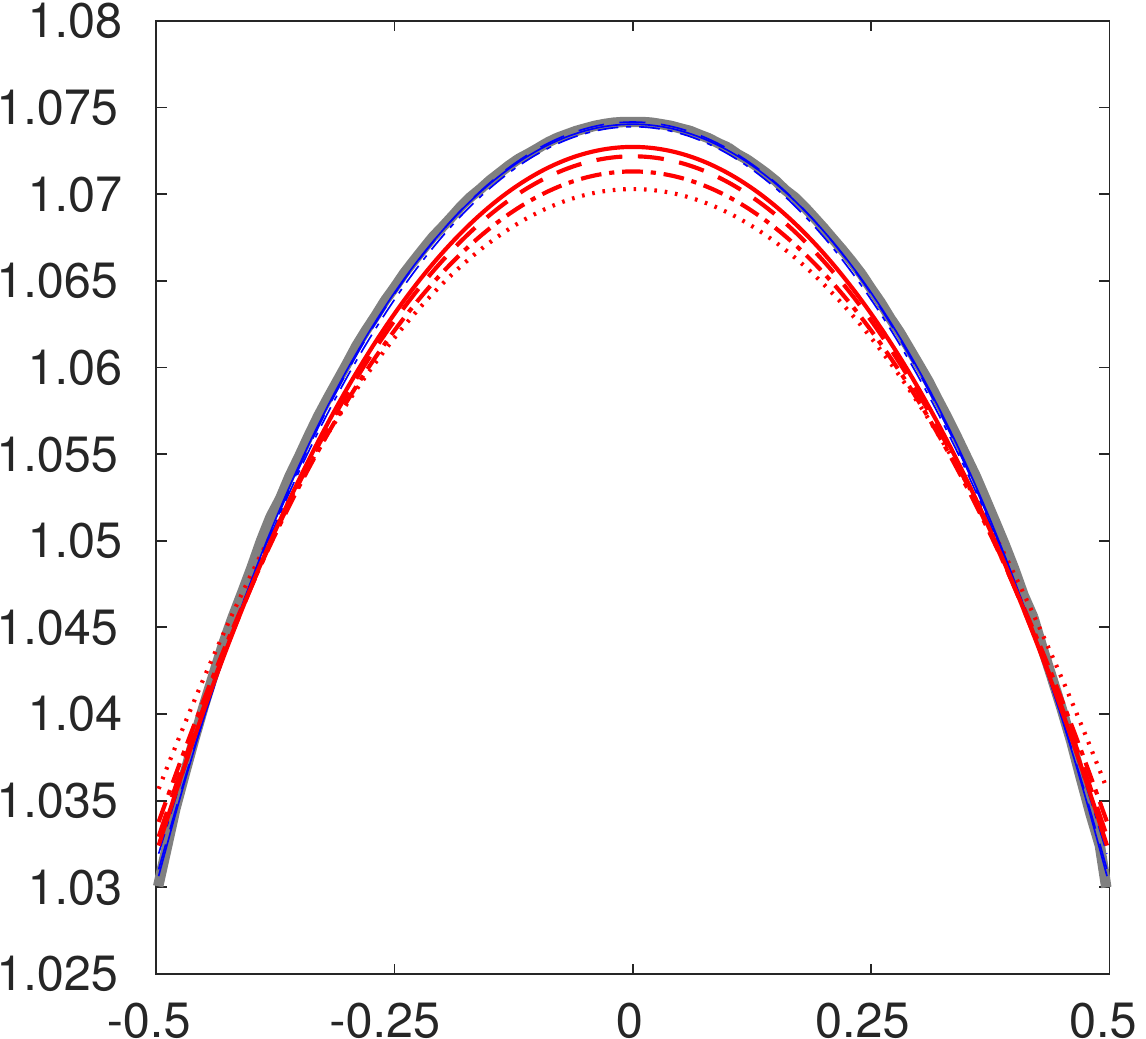}} \\
  \subfloat[Shear stress, $\sigma_{12}$]{\includegraphics[width=0.505\textwidth]{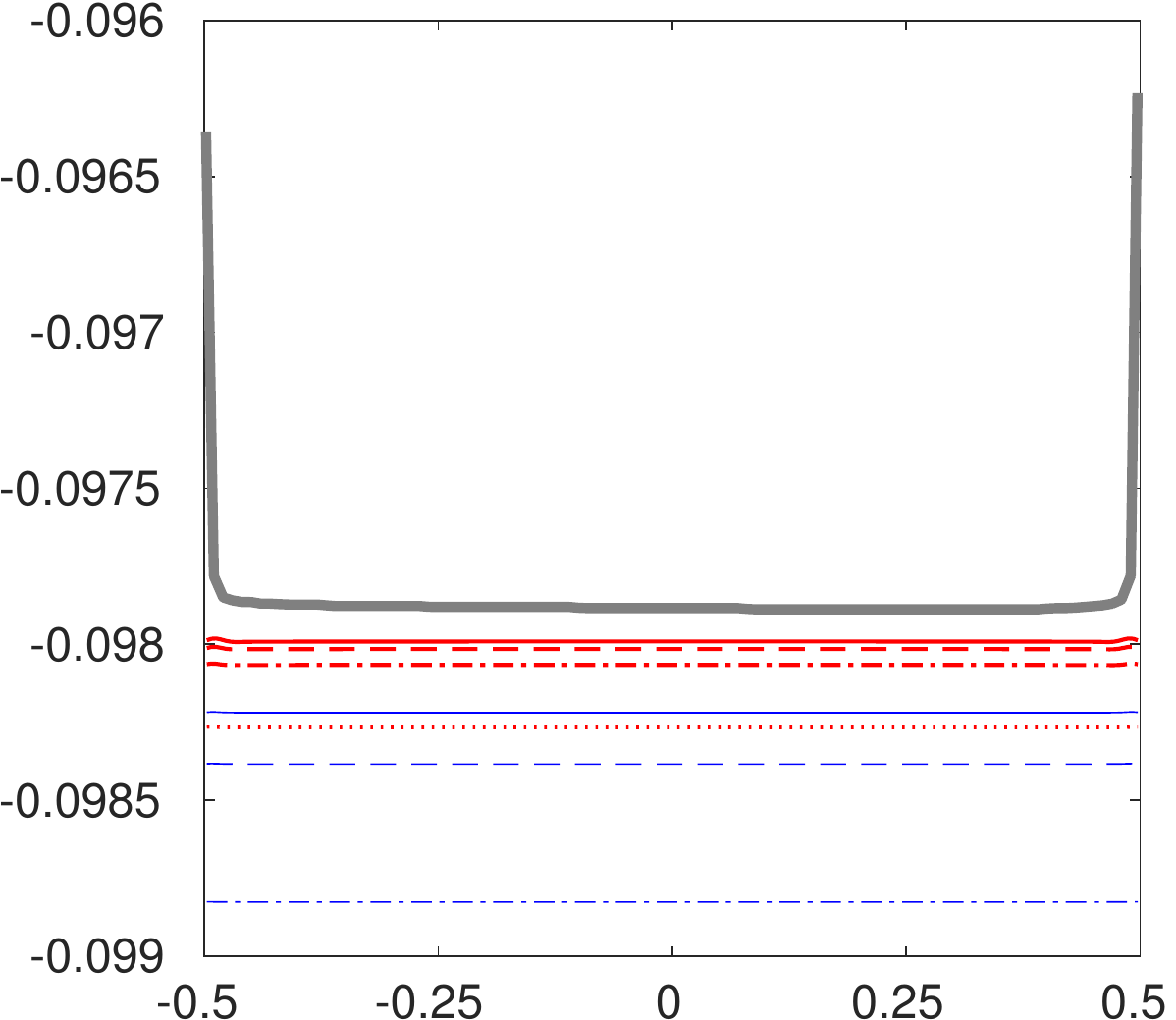}} \hfill
  \subfloat[Heat flux, $q_1$]{\includegraphics[width=0.48\textwidth]{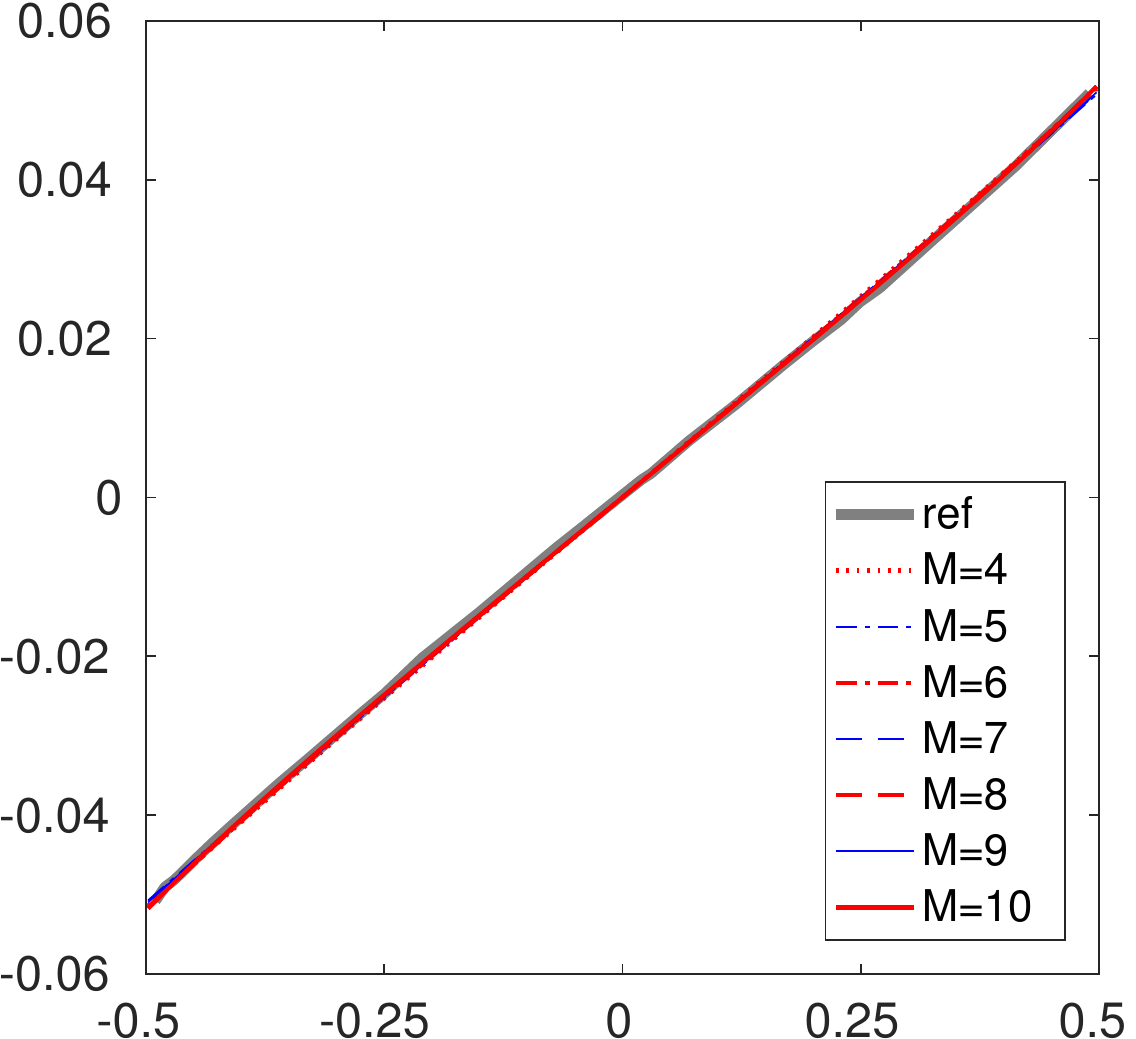}} 
\caption{Solution of the Couette flow for $\Kn=0.1199$ on the uniform grid with $N=200$.}
  \label{fig:couette-Kn01-uw300}
\end{figure}

\begin{figure}[!htb]
  \centering
  \subfloat[Density, $\rho$]{\includegraphics[width=0.49\textwidth]{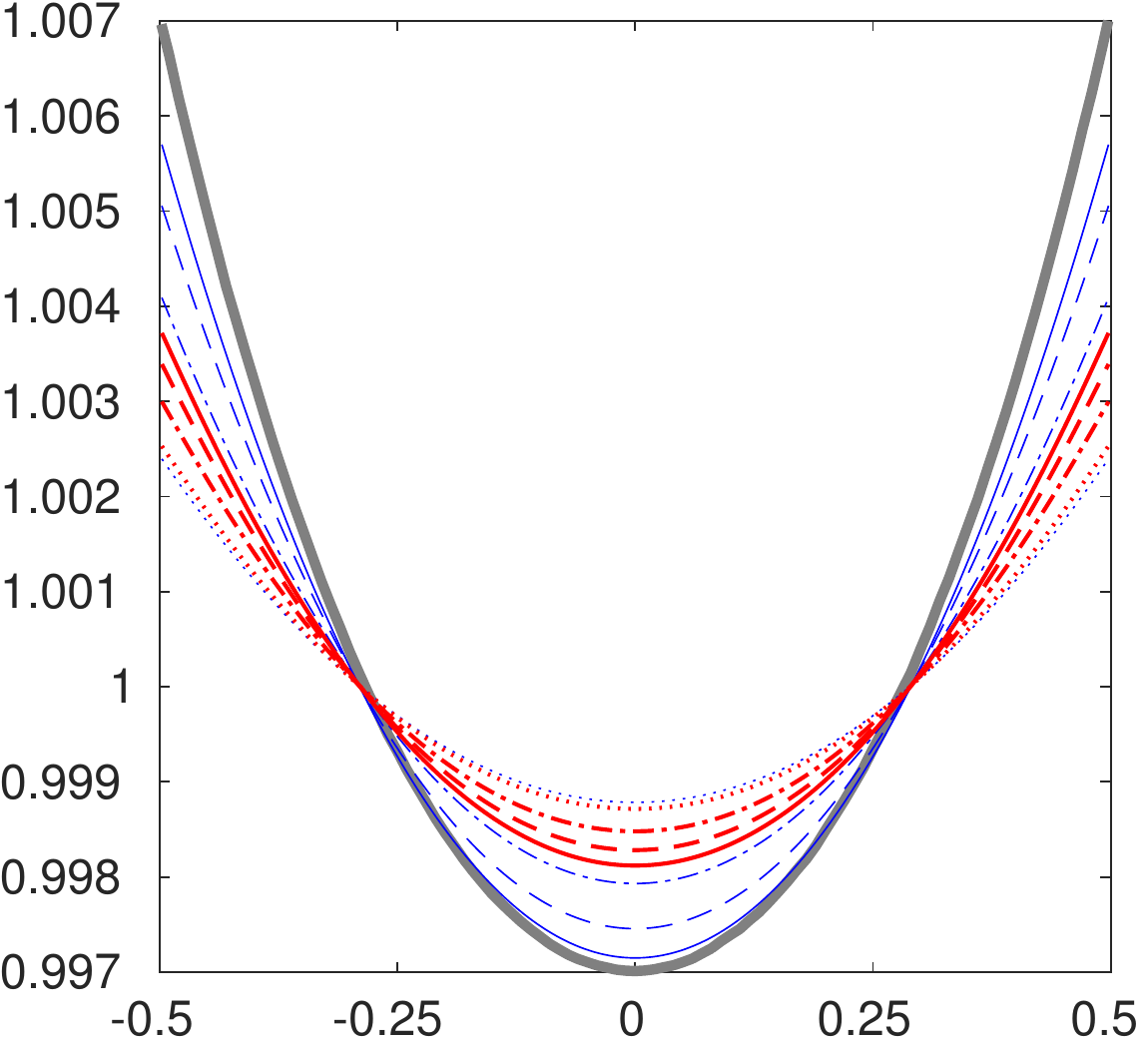}}\hfill
  \subfloat[Temperature, $\theta$]{\includegraphics[width=0.49\textwidth]{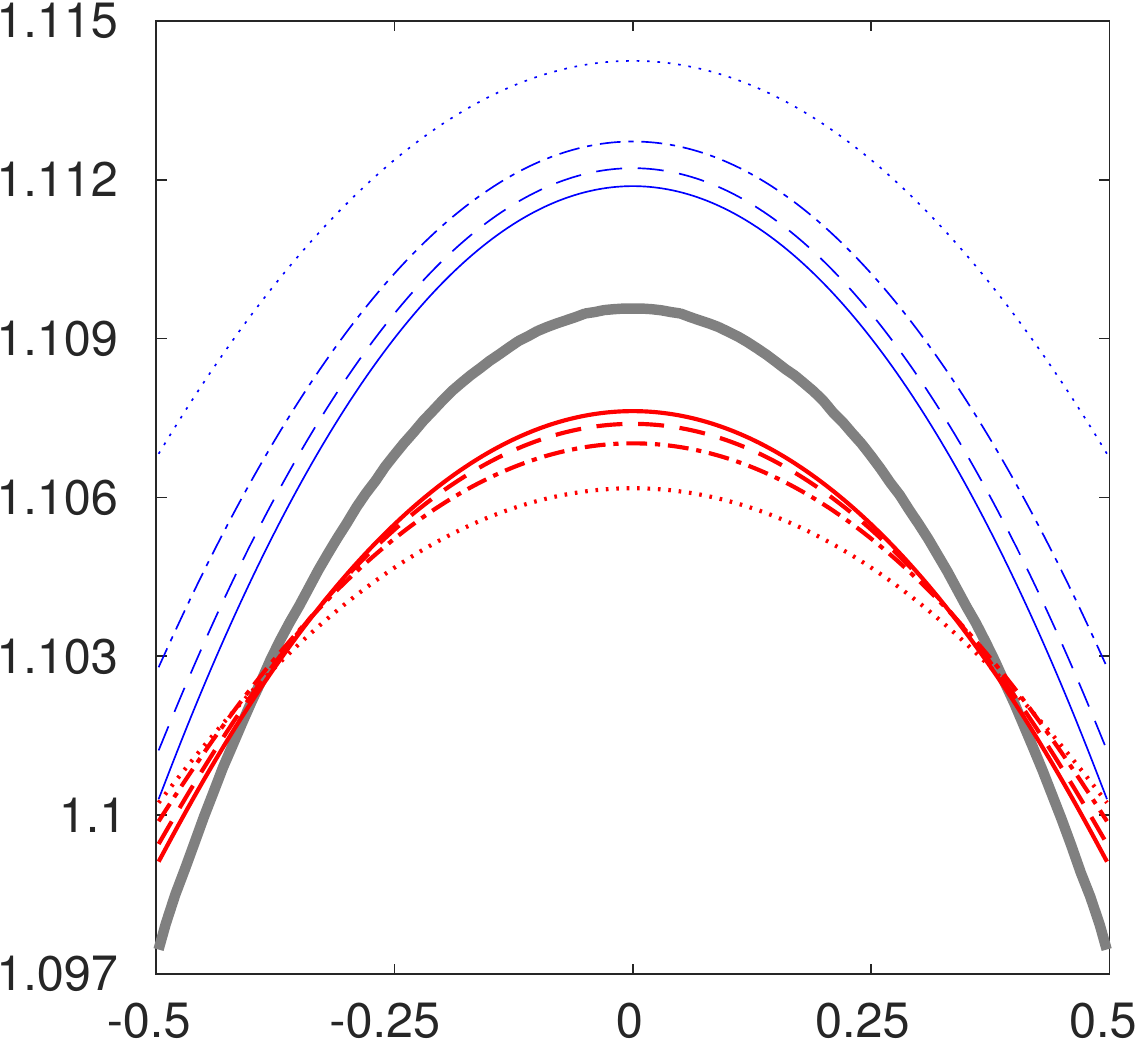}} \\
  \subfloat[Shear stress, $\sigma_{12}$]{\includegraphics[width=0.50\textwidth]{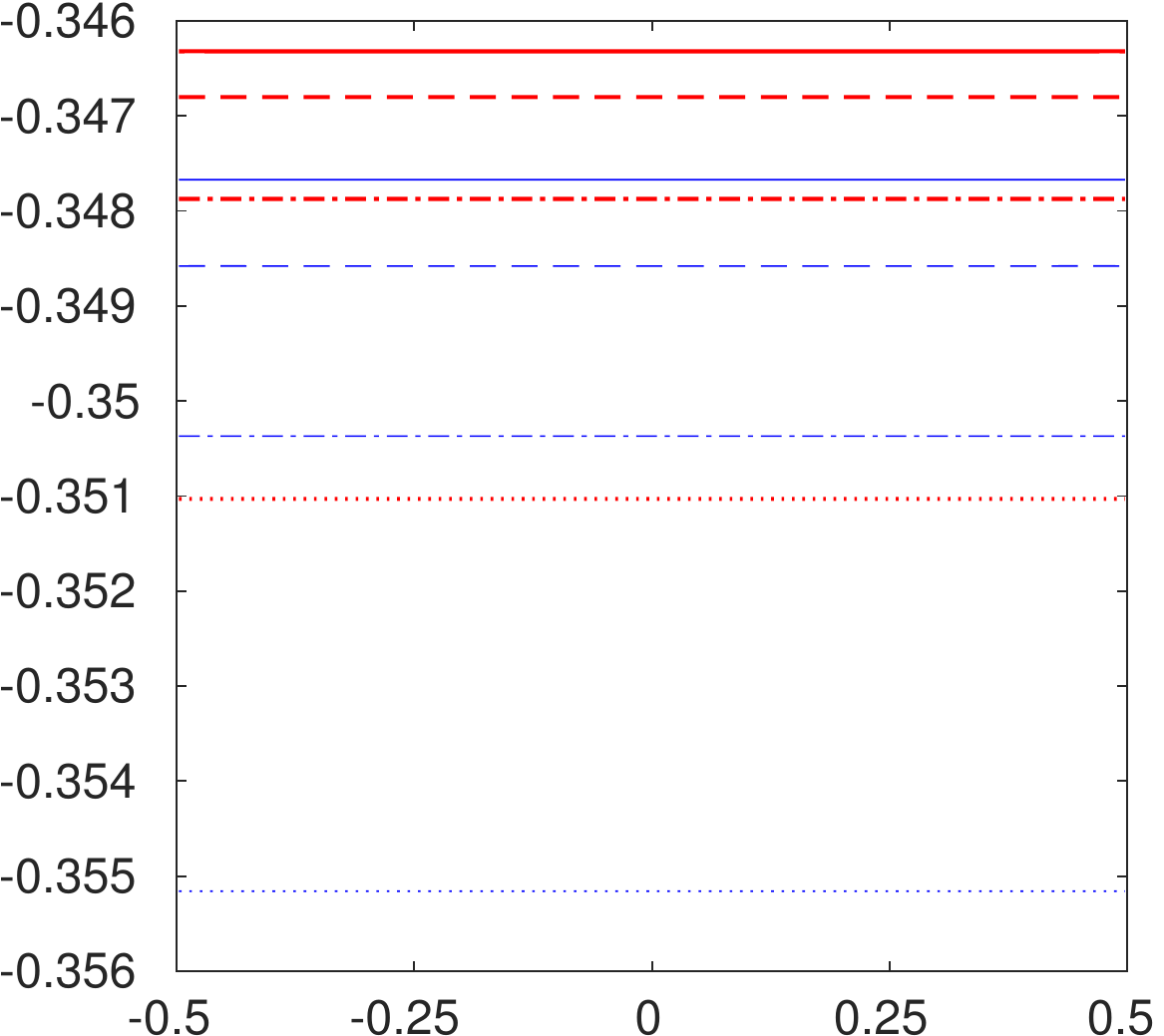}} \hfill
  \subfloat[Heat flux, $q_1$]{\includegraphics[width=0.49\textwidth]{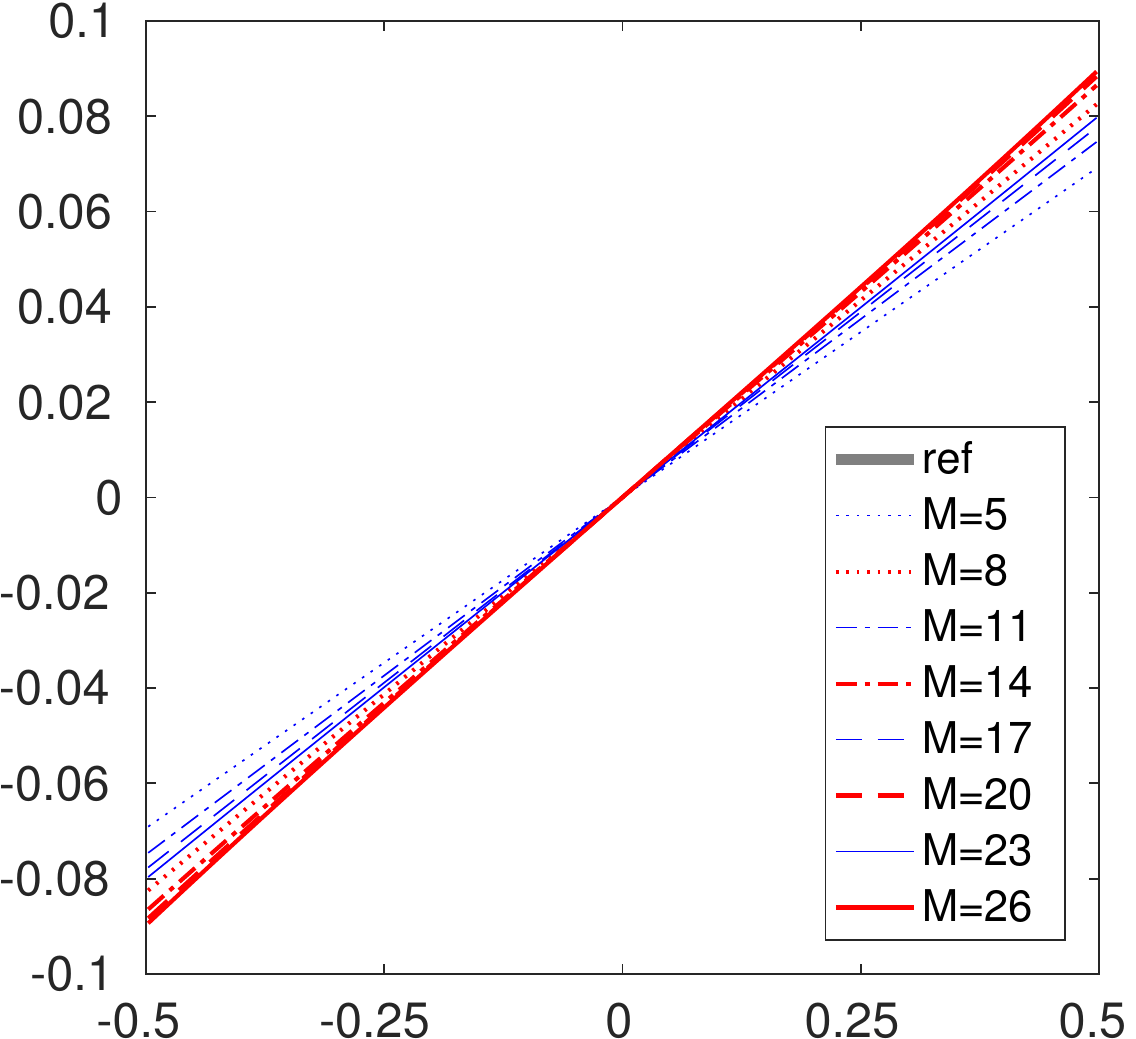}} 
\caption{Solution of the Couette flow for $\Kn=1.199$ on the uniform grid with $N=200$.}
  \label{fig:couette-Kn10-uw300}
\end{figure}

Numerical solutions of the moment models \eqref{eq:mnt-eqs} for
density $\rho$, temperature $\theta$, shear stress $\sigma_{12}$ and
heat flux $q_{1}$ on the uniform grid with $N=200$ cells are listed in
\figurename~\ref{fig:couette-Kn01-uw300} and
\ref{fig:couette-Kn10-uw300} for $\Kn = 0.1199$ and $1.199$
respectively. The solutions obtained by the discrete velocity method
\cite{Mieussens2004} are provided as a reference. We omit the
discussion on the accuracy and convergence of our results with respect
to $M$ here, since we are actually reproducing the results obtained in
\cite{Microflows1D}, where the validation of them has been
investigated in detail. We would just like to mention that the moment
model of order $M=10$ is sufficient to give satisfactory results for
$\Kn=0.1199$, while in the case $\Kn=1.199$ the moment model up to
order $M=23$ or $26$ is necessary. Moreover, a careful comparison
shows that the present second-order spatial discretization with
$N=200$ gives a slightly better results than those obtained in
\cite{Microflows1D, hu2014nmg} by the first-order spatial
discretization with $N=2048$, which indicates a remarkable improvement
in efficiency is obtained.

Now we turn to investigate the efficiency and behavior of the NMLM
solver proposed in previous sections. As we have done in
\cite{hu2016acceleration}, the NMLM solvers with various levels and
order reduction strategies for
the above Couette flow are performed on three uniform grids with
$N=100$, $200$ and $400$, respectively. Due to the similar features of
the NMLM solver with respect to $M$, only partial results are provided
here. 

In the case of $\Kn=0.1199$, the total number of iterations and the
elapsed CPU seconds, spent by the steady-state computation of the
solver, as well as the comparison to their counterparts of the single
level solver, are listed in
\tablename~\ref{tab:couette-Kn01-uw300:nmm-M45} for $M=4,5$ and in
\tablename~\ref{tab:couette-Kn01-uw300:nmm-M10-p1}-\ref{tab:couette-Kn01-uw300:nmm-M10-p2}
for $M=10$, where $K$ and $T$ denote respectively the total number of
iterations and the elapsed CPU seconds, and the corresponding
quantities of the single level solver are denoted by $K_{s}$ and
$T_{s}$. The convergence histories of the tests on the uniform grid
with $N=200$ are plotted in
\figurename~\ref{fig:couette-Kn01-uw300:nmm-res-history-N200-M4}-\ref{fig:couette-Kn01-uw300:nmm-res-history-N200-M10}
for $M=4, 5$ and $10$ respectively. It can be observed that the NMLM
solver, in comparison to the single level solver, could accelerate the
steady-state computation a lot for all tests. In more detail, the
total number of the NMLM iterations, for the same $M$ and the same
order reduction strategy, decreases as the total levels of the solver
increases, which indicates the convergence rate is
improved. Consequently, the elapsed CPU time is reduced 
as the total levels increases. For the NMLM solver with the same total
levels, the convergence rate of the order reduction strategy
$m_{\sss l-1} = \lceil m_{\sss l} / 2 \rceil$ is better than the
strategy $m_{\sss l-1} = m_{\sss l} - 2$, and the latter strategy is
better than the strategy $m_{\sss l-1} = m_{\sss l} - 1$. Apparently,
the computational cost of each NMLM iteration for these three order
reduction strategies is in the ascending sort, since smallest order is
employed in each level of the lower-order moment model correction for
the strategy $m_{\sss l-1} = \lceil m_{\sss l} / 2 \rceil$, while
largest order is used in each level for the strategy
$m_{\sss l-1} = m_{\sss l} - 1$. Therefore, among these three order
reduction strategies, the most efficient one becomes
$m_{\sss l-1} = \lceil m_{\sss l} / 2 \rceil$, the second one is
$m_{\sss l-1} = m_{\sss l} - 2$, and the third one is
$m_{\sss l-1} = m_{\sss l} - 1$. Although more total levels can be
applied for the strategies $m_{\sss l-1} = m_{\sss l} - 2$ and
$m_{\sss l-1} = m_{\sss l} - 1$, they might still not be as efficient
as the strategy $m_{\sss l-1} = \lceil m_{\sss l} / 2 \rceil$. Taking
the tests for $M=10$ as an example, the $3$-level NMLM solver with the
strategy $m_{\sss l-1} = \lceil m_{\sss l} / 2 \rceil$ is comparable
with $5$-level NMLM solver with the strategy
$m_{\sss l-1} = m_{\sss l} - 2$, and both of them are more efficient
than the $8$-level NMLM solver with the strategy
$m_{\sss l-1} = m_{\sss l} - 1$. At last, it can also be observed from
\tablename~\ref{tab:couette-Kn01-uw300:nmm-M45}-\ref{tab:couette-Kn01-uw300:nmm-M10-p2}
that the behavior of the NMLM solver with respect to the spatial grid
number $N$ is just like the one of the single level solver, that is,
the total number of NMLM iterations doubles and the elapsed CPU time
quadruples, as $N$ doubles.

As the Knudsen number increases to $\Kn=1.199$, the moment model with
a larger order needs to be considered. A partial numerical results can
be found in \tablename~\ref{tab:couette-Kn10-uw300:nmm-M23-p2} for
$M=23$ and in \tablename~\ref{tab:couette-Kn10-uw300:nmm-M26-p2} for
$M=26$, respectively. The corresponding convergence histories on the
uniform grid with $N=200$ are shown in
\figurename~\ref{fig:couette-Kn10-uw300:nmm-res-history-N200-M23}-\ref{fig:couette-Kn10-uw300:nmm-res-history-N200-M26}. Like
the case of $\Kn=0.1199$, it can be observed that the NMLM solver
behaves similarly to the corresponding single level solver, as the
grid number $N$ increases. The convergence rate of the NMLM solver
with the same order reduction strategies is also improved as the total
levels increases. The elapsed CPU time is consequently reduced except
for the tests with $M=26$ and $m_{\sss l-1} = m_{\sss l} - 2$, which
is acceptable by noting that the convergence rate is improved a
little, and the computational cost of lower-order moment model
correction at each level could not be neglected, since the order
sequence $26, 24, 22, \ldots$ is adopted. Moreover, oscillation of the
residual at the beginning iterations and degeneracy of the convergence
rate are observed for the single level solver, i.e., Heun's
method. This makes the residual oscillate more wildly and the
convergence rate also be degenerated for the NMLM solver, especially
for the solver with the strategy
$m_{\sss l-1} = \lceil m_{\sss l} / 2 \rceil$, for which the
convergence rate, in contrast to the case of $\Kn=0.1199$, is now
worse than the strategy $m_{\sss l-1} = m_{\sss l} - 2$. Nevertheless,
due to the great reduction of the computational cost at each NMLM
iteration, the strategy $m_{\sss l-1} = \lceil m_{\sss l} / 2 \rceil$
is finally more efficient than the strategy
$m_{\sss l-1} = m_{\sss l} - 2$. As can be seen, more than $50\%$ of
the total computational cost, compared with the single level solver,
is saved by the $4$-level NMLM solver with the strategy
$m_{\sss l-1} = \lceil m_{\sss l} / 2 \rceil$ in all tests for $M=23$
and $26$. In addition, to seek the balance between the convergence
rate and the computational cost of each NMLM iteration, a new order
reduction strategy, namely, $m_{\sss l-1} = m_{\sss l} - 4$, is tested
for $M=23$. As shown in
\figurename~\ref{fig:couette-Kn10-uw300:nmm-res-history-N200-M23}, it
is found that the convergence rate of this strategy is better than the
strategy $m_{\sss l-1} = \lceil m_{\sss l} / 2 \rceil$, and the
elapsed CPU time of the $6$-level NMLM solver with the former strategy
is a slightly less than the $4$-level NMLM solver with the latter
strategy.

\begin{figure}[!htb]
  \centering
  \subfloat[$\Kn=0.1199$]{\includegraphics[width=0.49\textwidth]{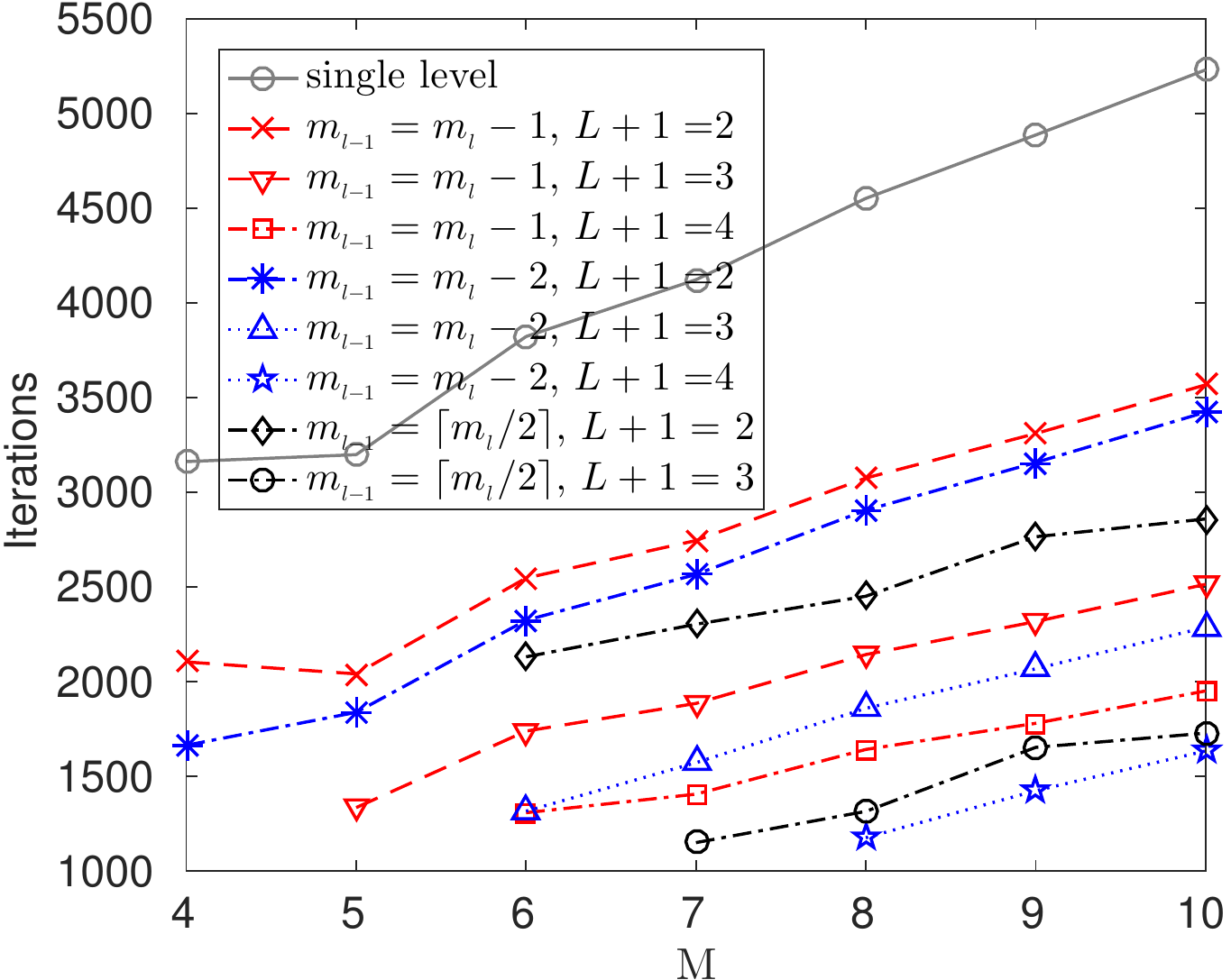}}\hfill
  \subfloat[$\Kn=1.199$]{\includegraphics[width=0.49\textwidth]{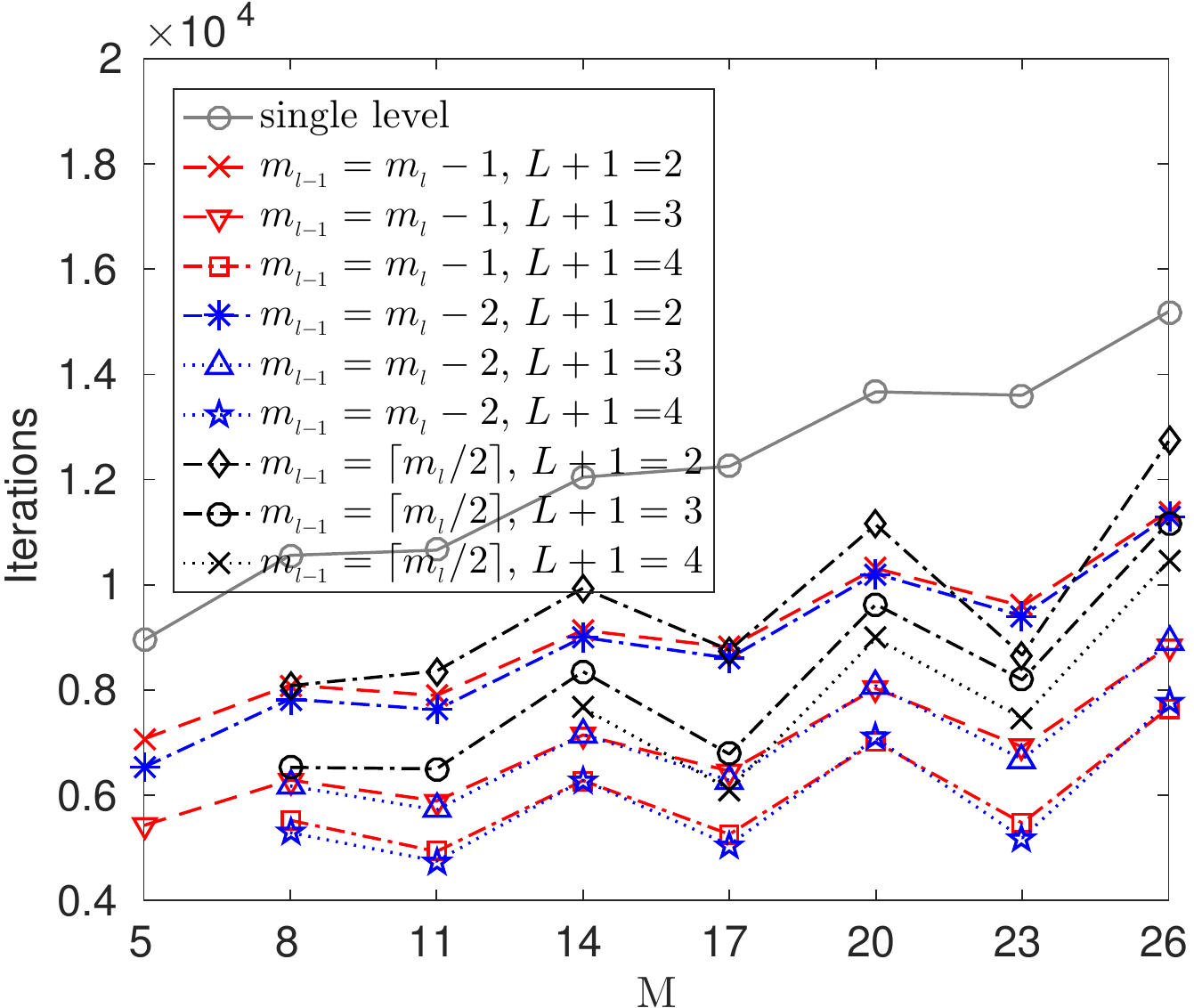}}
  \caption{Total number of iterations in terms of $M$ of the NMLM
    solver for the Couette flow on the uniform grid of $N=200$. The
    total number of iterations of the single level solver is rescaled
    by a factor of $6$.}
  \label{fig:couette:nmm-iters-orders}
\end{figure}

Finally, the behavior of the total number of iterations with respect
to the order of the moment model $M$ is investigated. The results are
shown in \figurename~\ref{fig:couette:nmm-iters-orders}. It can be
seen that in the case of $\Kn=0.1199$, the total number of iterations
increases almost linearly with the same ratio for all tests, as $M$
increases. While in the case of $\Kn=1.199$, sawtooth polylines are
observed for all tests. To be specific, the single level solver for
odd $M$ performs much better than the solver for successor even $M$,
and the growth rate of the total number of iterations with respect to
odd or even $M$ is nearly the same. This shows better performance of
Heun's method than the SGS-Richardson iteration, in comparison to the
results presented in \cite{hu2016acceleration}. For the multi-level
NMLM solver, the different performance for odd or even $M$ becomes
more obvious, especially for the strategy
$m_{\sss l-1} = \lceil m_{\sss l} / 2 \rceil$. The underlying reason
remains to be further studied. However, we can observe that the growth
rate of the total number of iterations for the multi-level NMLM solver
with respect to even $M$ is greater than the corresponding growth rate
with respect to odd $M$, but still be almost not greater than the
growth rate of the single level solver. As a result, the NMLM solver
becomes more efficient for the moment model of odd $M$ than that of
even $M$.

\begin{table}[!ht]
  \centering\footnotesize
  \begin{tabular}{c|c||c|c|c||c|cc|c}
    \hline\hline
    \multicolumn{2}{c||}{} & \multicolumn{3}{c||}{$M=4$} & \multicolumn{4}{c}{$M=5$}\\
    \hline
    \multicolumn{2}{c||}{} & & \multicolumn{1}{c|}{$m_{\sss l-1} = m_{\sss l} - 1$} & \multicolumn{1}{c||}{$m_{\sss l-1} = m_{\sss l} - 2$} & & \multicolumn{2}{c|}{$m_{\sss l-1} = m_{\sss l} - 1$} & \multicolumn{1}{c}{$m_{\sss l-1} = m_{\sss l} - 2$}\\ 
    \hline
    \multicolumn{2}{c||}{$L+1$} & 1 & 2 & 2 & 1 & 2 & 3 & 2 \\
    \hline
    \multirow{4}{*}{\begin{sideways}$N=100$\end{sideways}} 
    & $K$          &        9484&        1052&         834&        9595&        1021&         668&         919\\
    & $T$          &      74.711&      51.155&      33.326&     124.430&      78.879&      62.440&      60.092\\
    & $K_s/K$      &       1.000&       9.015&      11.372&       1.000&       9.398&      14.364&      10.441\\
    & $T_s/T$      &       1.000&       1.460&       2.242&       1.000&       1.577&       1.993&       2.071\\
    \hline
    \multirow{4}{*}{\begin{sideways}$N=200$\end{sideways}} 
    & $K$          &       18971&        2104&        1664&       19191&        2041&        1335&        1838\\
    & $T$          &     317.736&     204.048&     132.755&     504.057&     315.070&     248.836&     238.401\\
    & $K_s/K$      &       1.000&       9.017&      11.401&       1.000&       9.403&      14.375&      10.441\\
    & $T_s/T$      &       1.000&       1.557&       2.393&       1.000&       1.600&       2.026&       2.114\\
    \hline
    \multirow{4}{*}{\begin{sideways}$N=400$\end{sideways}} 
    & $K$          &       37944&        4208&        3345&       38382&        4082&        2668&        3675\\
    & $T$          &    1278.954&     807.967&     538.051&    2028.431&    1260.130&     982.810&     955.244\\
    & $K_s/K$      &       1.000&       9.017&      11.343&       1.000&       9.403&      14.386&      10.444\\
    & $T_s/T$      &       1.000&       1.583&       2.377&       1.000&       1.610&       2.064&       2.123\\
    \hline\hline
  \end{tabular}
  \caption{Performance of the NMLM solver for the Couette flow with $\Kn=0.1199$ and $M=4,5$.}
  \label{tab:couette-Kn01-uw300:nmm-M45}
\end{table}

\begin{table}[!ht]
  \centering\footnotesize
  \begin{tabular}{c|c||ccccccc}
    \hline\hline
    \multicolumn{2}{c||}{} & \multicolumn{6}{c}{$m_{\sss l-1} = m_{\sss l} - 1$}\\
    \hline
    \multicolumn{2}{c||}{$L+1$} & 2 & 3 & 4 & 5 & 6 & 7 & 8 \\
    \hline
    \multirow{4}{*}{\begin{sideways}$N=100$\end{sideways}} 
    & $K$          &        1784&        1258&         976&         799&         675&         581&         506\\
    & $T$          &     883.677&     781.440&     704.394&     623.006&     553.480&     492.187&     442.410\\
    & $K_s/K$      &       8.796&      12.474&      16.078&      19.640&      23.247&      27.009&      31.012\\
    & $T_s/T$      &       1.336&       1.511&       1.676&       1.895&       2.133&       2.399&       2.669\\
    \hline
    \multirow{4}{*}{\begin{sideways}$N=200$\end{sideways}} 
    & $K$          &        3567&        2514&        1951&        1595&        1348&        1159&        1010\\
    & $T$          &    3546.025&    3164.404&    2791.258&    2521.665&    2248.097&    1988.022&    1734.307\\
    & $K_s/K$      &       8.798&      12.483&      16.085&      19.675&      23.280&      27.077&      31.071\\
    & $T_s/T$      &       1.349&       1.511&       1.713&       1.897&       2.127&       2.406&       2.758\\
    \hline
    \multirow{4}{*}{\begin{sideways}$N=400$\end{sideways}} 
    & $K$          &        7133&        5027&        3900&        3189&        2693&        2316&        2017\\
    & $T$          &   13688.115&   12610.887&   11097.524&   10001.285&    8973.447&    7835.308&    7006.392\\
    & $K_s/K$      &       8.799&      12.485&      16.093&      19.680&      23.305&      27.099&      31.116\\
    & $T_s/T$      &       1.332&       1.446&       1.643&       1.824&       2.032&       2.328&       2.603\\
    \hline\hline
  \end{tabular}
  \caption{Performance of the NMLM solver for the Couette flow with $\Kn=0.1199$ and $M=10$ (part I).}
  \label{tab:couette-Kn01-uw300:nmm-M10-p1}
\end{table}

\begin{table}[!ht]
  \centering\footnotesize
  \begin{tabular}{c|c||cccc|cc|c}
    \hline\hline
    \multicolumn{2}{c||}{} & \multicolumn{4}{c|}{$m_{\sss l-1} = m_{\sss l} - 2$} & \multicolumn{2}{c|}{$m_{\sss l-1} = \lceil m_{\sss l} / 2 \rceil$} & \\
    \hline
    \multicolumn{2}{c||}{$L+1$} & 2 & 3 & 4 & 5 & 2 & 3 & 1\\
    \hline
    \multirow{4}{*}{\begin{sideways}$N=100$\end{sideways}} 
    & $K$          &        1712&        1144&         818&         560&        1430&         864&       15692\\
    & $T$          &     739.831&     571.578&     421.999&     294.624&     478.640&     297.862&    1180.762\\
    & $K_s/K$      &       9.166&      13.717&      19.183&      28.021&      10.973&      18.162&       1.000\\
    & $T_s/T$      &       1.596&       2.066&       2.798&       4.008&       2.467&       3.964&       1.000\\
    \hline
    \multirow{4}{*}{\begin{sideways}$N=200$\end{sideways}} 
    & $K$          &        3423&        2287&        1634&        1118&        2859&        1728&       31382\\
    & $T$          &    2940.028&    2250.739&    1694.924&    1165.652&    1923.646&    1155.322&    4782.708\\
    & $K_s/K$      &       9.168&      13.722&      19.206&      28.070&      10.977&      18.161&       1.000\\
    & $T_s/T$      &       1.627&       2.125&       2.822&       4.103&       2.486&       4.140&       1.000\\
    \hline
    \multirow{4}{*}{\begin{sideways}$N=400$\end{sideways}} 
    & $K$          &        6845&        4572&        3265&        2234&        5716&        3454&       62761\\
    & $T$          &   11978.715&    8947.557&    6792.666&    4643.324&    7329.159&    4726.841&   18238.107\\
    & $K_s/K$      &       9.169&      13.727&      19.222&      28.094&      10.980&      18.171&       1.000\\
    & $T_s/T$      &       1.523&       2.038&       2.685&       3.928&       2.488&       3.858&       1.000\\
    \hline\hline
  \end{tabular}
  \caption{Performance of the NMLM solver for the Couette flow with $\Kn=0.1199$ and $M=10$ (part II).}
  \label{tab:couette-Kn01-uw300:nmm-M10-p2}
\end{table}

\begin{table}[!ht]
  \hspace*{-1.em}
  \centering\scriptsize
  \begin{tabular}{c|c||ccccc|ccc}
    \hline\hline
    \multicolumn{2}{c||}{} & \multicolumn{5}{c|}{$m_{\sss l-1} = m_{\sss l} - 2$} & \multicolumn{3}{c}{$m_{\sss l-1} = \lceil m_{\sss l} / 2 \rceil$} \\
    \hline
    \multicolumn{2}{c||}{$L+1$} & 4 & 5 & 6 & 7 & 8 & 2 & 3 & 4 \\
    \hline
    \multirow{4}{*}{\begin{sideways}$N=100$\end{sideways}} 
    & $K$          &        2580&        2194&        2033&        1910&        1802&        4198&        4111&        3737\\
    & $T$          &   21771.636&   20636.803&   20038.953&   18945.630&   18761.918&   16608.830&   16448.863&   14880.760\\
    & $K_s/K$      &      15.820&      18.603&      20.076&      21.369&      22.650&       9.722&       9.928&      10.922\\
    & $T_s/T$      &       1.776&       1.874&       1.929&       2.041&       2.061&       2.328&       2.351&       2.598\\
    \hline
    \multirow{4}{*}{\begin{sideways}$N=200$\end{sideways}} 
    & $K$          &        5188&        4427&        4096&        3853&        3644&        8665&        8195&        7465\\
    & $T$          &   88518.030&   82354.633&   80539.022&   77774.608&   74883.559&   68274.318&   65399.218&   59797.985\\
    & $K_s/K$      &      15.733&      18.437&      19.927&      21.184&      22.399&       9.420&       9.960&      10.934\\
    & $T_s/T$      &       1.746&       1.877&       1.919&       1.987&       2.064&       2.264&       2.363&       2.584\\
    \hline
    \multirow{4}{*}{\begin{sideways}$N=400$\end{sideways}} 
    & $K$          &       10396&        8875&        8203&        7715&        7298&       17519&       16362&       14896\\
    & $T$          &  351768.755&  332364.873&  323203.956&  313285.091&  300010.067&  274114.156&  258269.847&  236863.512\\
    & $K_s/K$      &      15.702&      18.393&      19.900&      21.159&      22.368&       9.318&       9.977&      10.959\\
    & $T_s/T$      &       1.743&       1.845&       1.897&       1.957&       2.044&       2.237&       2.374&       2.588\\
    \hline\hline
  \end{tabular}
  \caption{Performance of the NMLM solver for the Couette flow with $\Kn=1.199$ and $M=23$.}
  \label{tab:couette-Kn10-uw300:nmm-M23-p2}
\end{table}

\begin{table}[!ht]
\hspace*{-4em}
  \centering\scriptsize
  \begin{tabular}{c|c||ccccc|cccc}
    \hline\hline
    \multicolumn{2}{c||}{} & \multicolumn{5}{c|}{$m_{\sss l-1} = m_{\sss l} - 2$} & \multicolumn{4}{c}{$m_{\sss l-1} = \lceil m_{\sss l} / 2 \rceil$} \\
    \hline
    \multicolumn{2}{c||}{$L+1$} & 4 & 5 & 6 & 7 & 8 & 2 & 3 & 4 & 5\\
    \hline
    \multirow{4}{*}{\begin{sideways}$N=100$\end{sideways}} 
    & $K$          &        3882&        3567&        3372&        3239&        3143&        6368&        5594&        5234&        5042\\
    & $T$          &   48452.216&   49497.714&   49879.413&   50491.512&   50094.156&   36028.855&   31647.628&   29804.057&   28720.545\\
    & $K_s/K$      &      11.784&      12.825&      13.566&      14.123&      14.555&       7.184&       8.178&       8.740&       9.073\\
    & $T_s/T$      &       1.267&       1.240&       1.231&       1.216&       1.225&       1.704&       1.940&       2.060&       2.137\\
    \hline
    \multirow{4}{*}{\begin{sideways}$N=200$\end{sideways}} 
    & $K$          &        7768&        7143&        6753&        6488&        6295&       12731&       11183&       10457&       10064\\
    & $T$          &  193489.560&  199169.100&  199861.340&  202753.586&  199230.828&  142907.221&  126167.456&  118536.002&  114107.992\\
    & $K_s/K$      &      11.736&      12.763&      13.500&      14.051&      14.482&       7.161&       8.152&       8.718&       9.058\\
    & $T_s/T$      &       1.268&       1.232&       1.227&       1.210&       1.231&       1.717&       1.944&       2.069&       2.150\\
    \hline
    \multirow{4}{*}{\begin{sideways}$N=400$\end{sideways}} 
    & $K$          &       15538&       14293&       13516&       12987&       12600&       25447&       22360&       20902&       20107\\
    & $T$          &  772099.850&  792816.382&  801340.202&  805215.081&  797777.655&  573505.900&  506611.853&  473278.177&  458730.507\\
    & $K_s/K$      &      11.735&      12.757&      13.491&      14.040&      14.472&       7.166&       8.155&       8.724&       9.069\\
    & $T_s/T$      &       1.296&       1.262&       1.248&       1.242&       1.254&       1.744&       1.975&       2.114&       2.181\\
    \hline\hline
  \end{tabular}
  \caption{Performance of the NMLM solver for the Couette flow with $\Kn=1.199$ and $M=26$.}
  \label{tab:couette-Kn10-uw300:nmm-M26-p2}
\end{table}

\newcommand\drawCouetteNMMHistory[6]{\begin{figure}[!htb]
  \centering
  {\includegraphics[width=0.5\textwidth]{couette_#1_#3_nmm_res_iters_N#5_M#6.pdf}}\hfill
  {\includegraphics[width=0.5\textwidth]{couette_#1_#3_nmm_res_cputime_N#5_M#6.pdf}}
  \caption{Convergence history of the NMLM solver for the Couette flow
    with $\Kn=#2$ and $M=#6$ on the uniform grid of $N=#5$.}
  \label{fig:couette-#1-#3:nmm-res-history-N#5-M#6}
\end{figure}
}
\drawCouetteNMMHistory{Kn01}{0.1199}{uw300}{}{200}{4}
\drawCouetteNMMHistory{Kn01}{0.1199}{uw300}{}{200}{5}
\drawCouetteNMMHistory{Kn01}{0.1199}{uw300}{}{200}{10}
\drawCouetteNMMHistory{Kn10}{1.199}{uw300}{}{200}{23}
\drawCouetteNMMHistory{Kn10}{1.199}{uw300}{}{200}{26}

\subsection{The force driven Poiseuille flow}
\label{sec:num-ex-poiseuille}
The second example is the force driven Poiseuille flow which has been
investigated in the literature, see e.g. \cite{Li, Xu2007, Garcia}.
The gas lies between two infinite parallel plates which are stationary
and have the same temperature of $\theta^{W}=1$. It is driven by an
external constant force and has a steady state as time goes. In
our simulation, the distance of the two plates is assumed to be
$L_{D}=1$, and the acceleration due to the external force is set to be
$\bF=(0, 0.2555,0)^T$. The collision frequency for the variable hard
sphere model, that is,
\begin{align}
  \label{eq:vhs-nu}
  \nu = \sqrt{\frac{2}{\pi}} \frac{(5-2w)(7-2w)\Pr}{15 \Kn} \rho \theta^{1-w},
\end{align}
with the viscosity index $w=0.5$ and the Knudsen number $\Kn=0.1$ is
adopted.

\begin{figure}[!htb]
  \centering
  \subfloat[Density, $\rho$]{\includegraphics[width=0.49\textwidth]{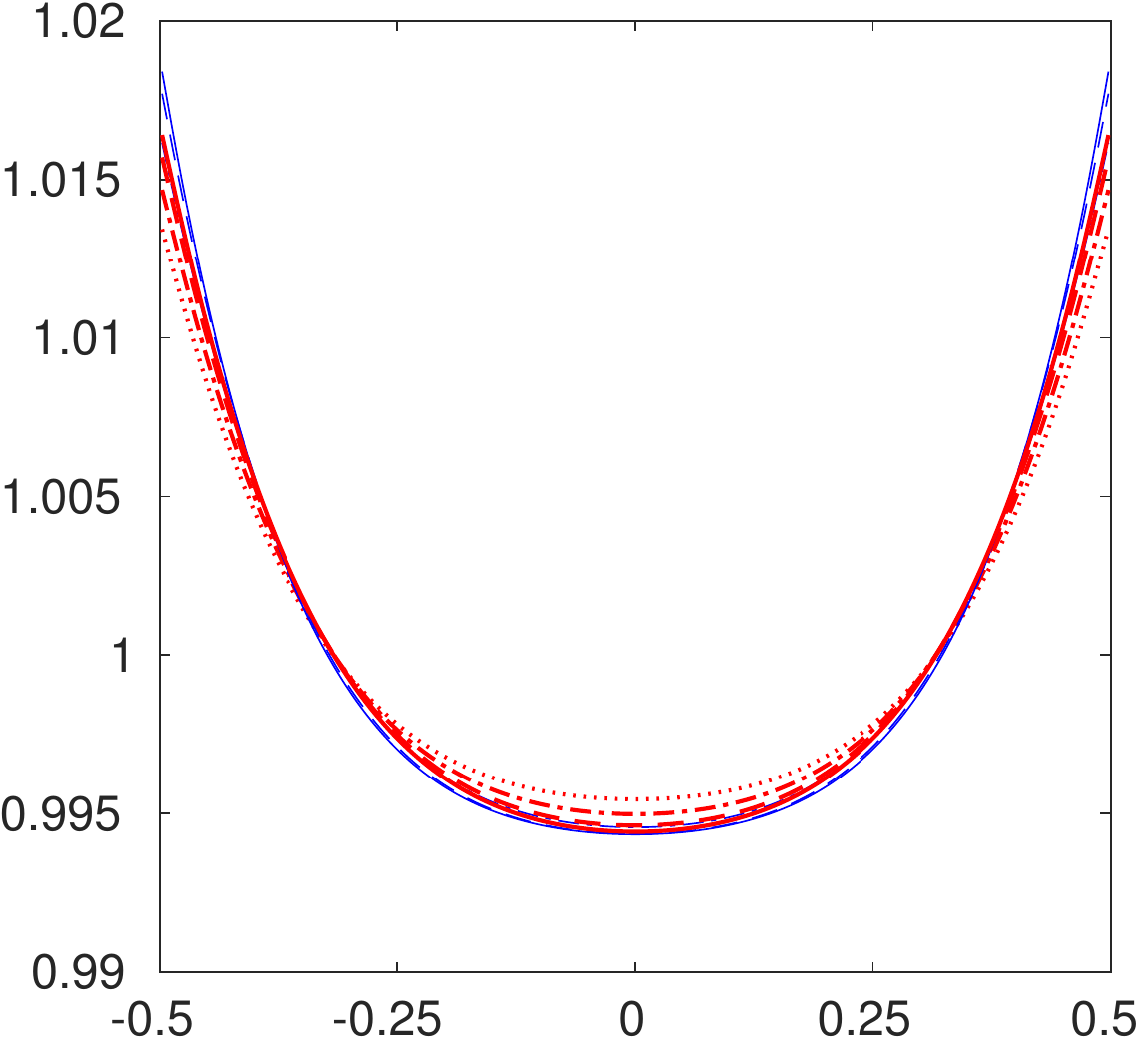}}\hfill
  \subfloat[Temperature, $\theta$]{\includegraphics[width=0.49\textwidth]{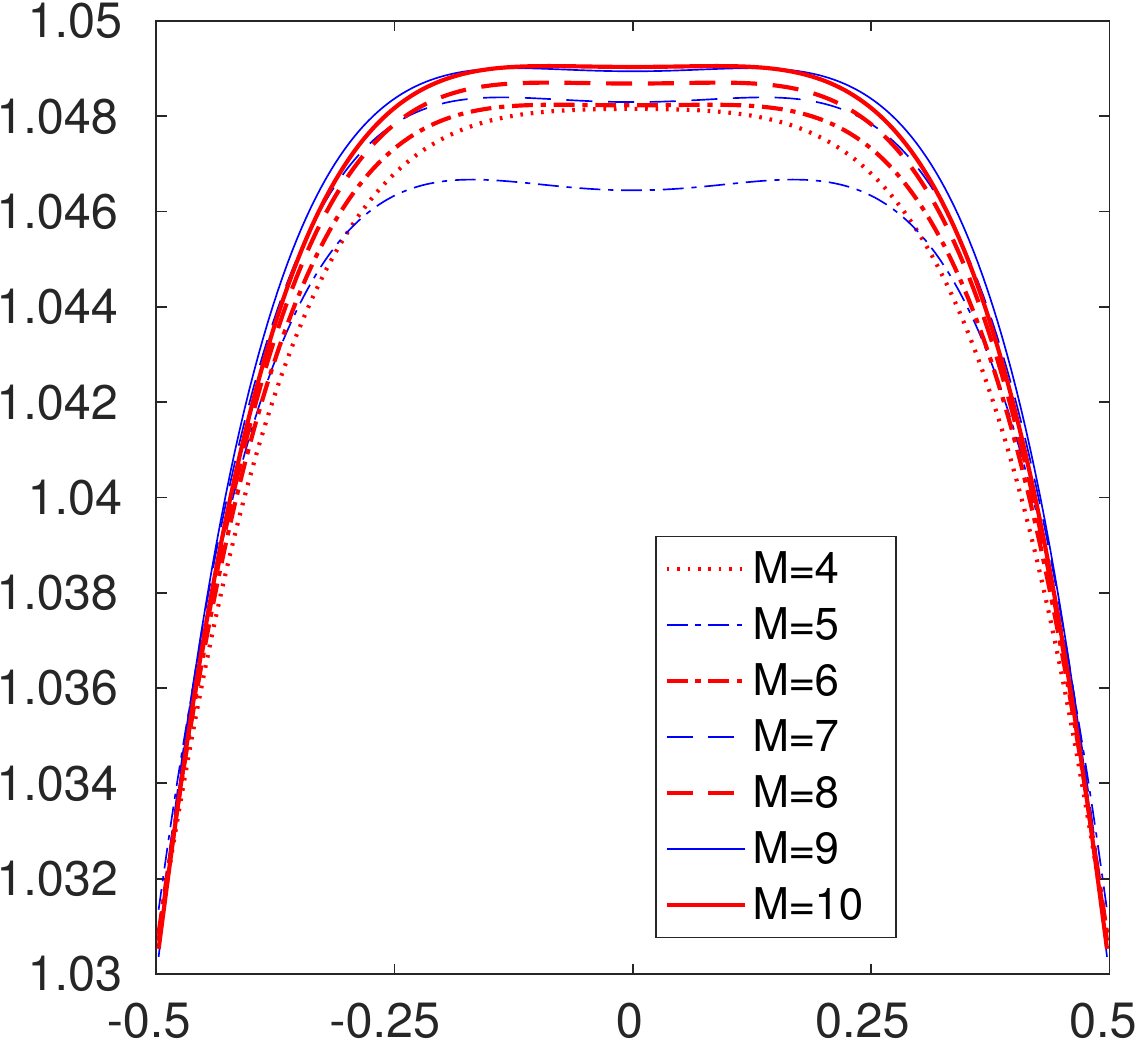}} \\
  \subfloat[Normal stress, $\sigma_{11}$]{\includegraphics[width=0.475\textwidth]{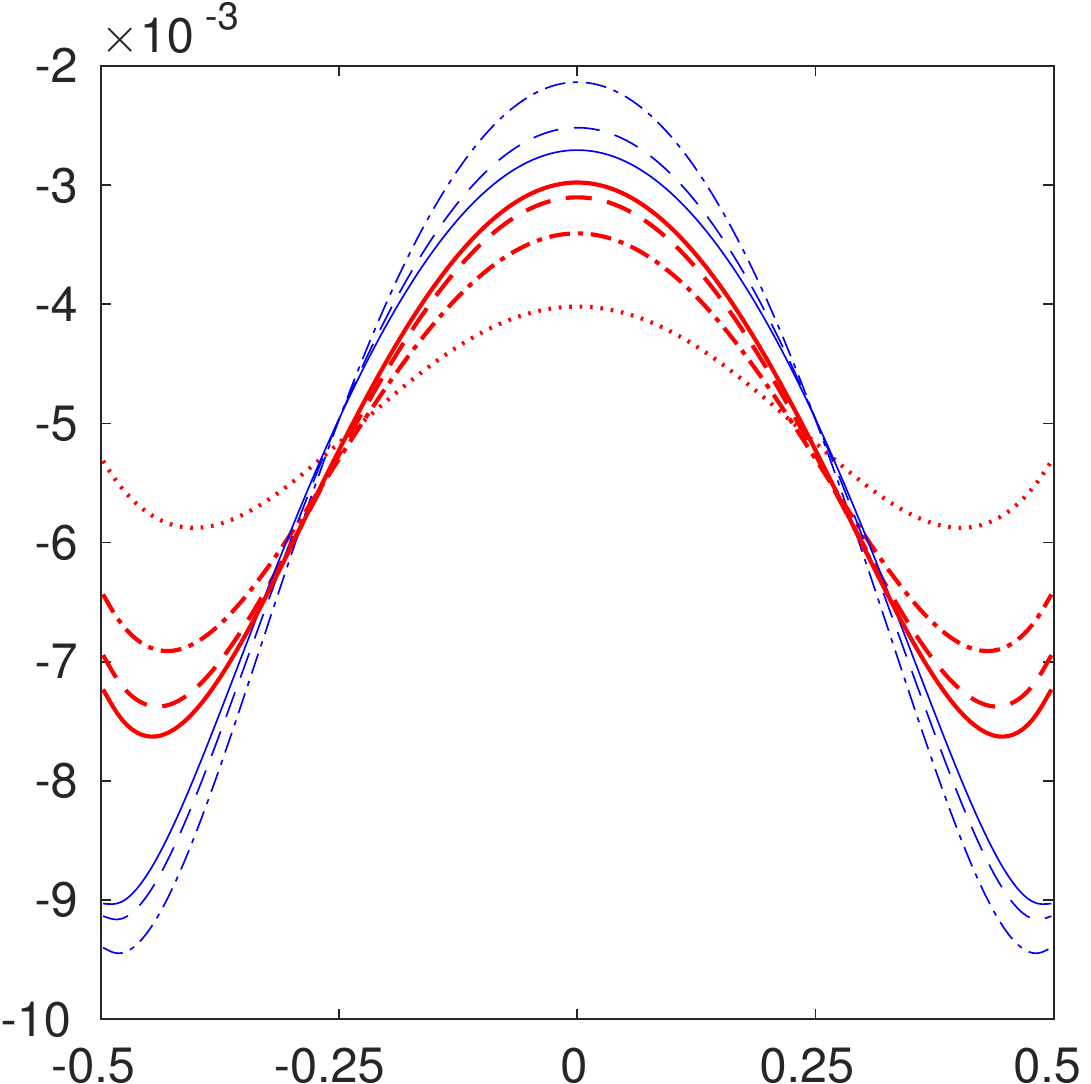}} \hfill
  \subfloat[Heat flux, $q_2$]{\includegraphics[width=0.51\textwidth]{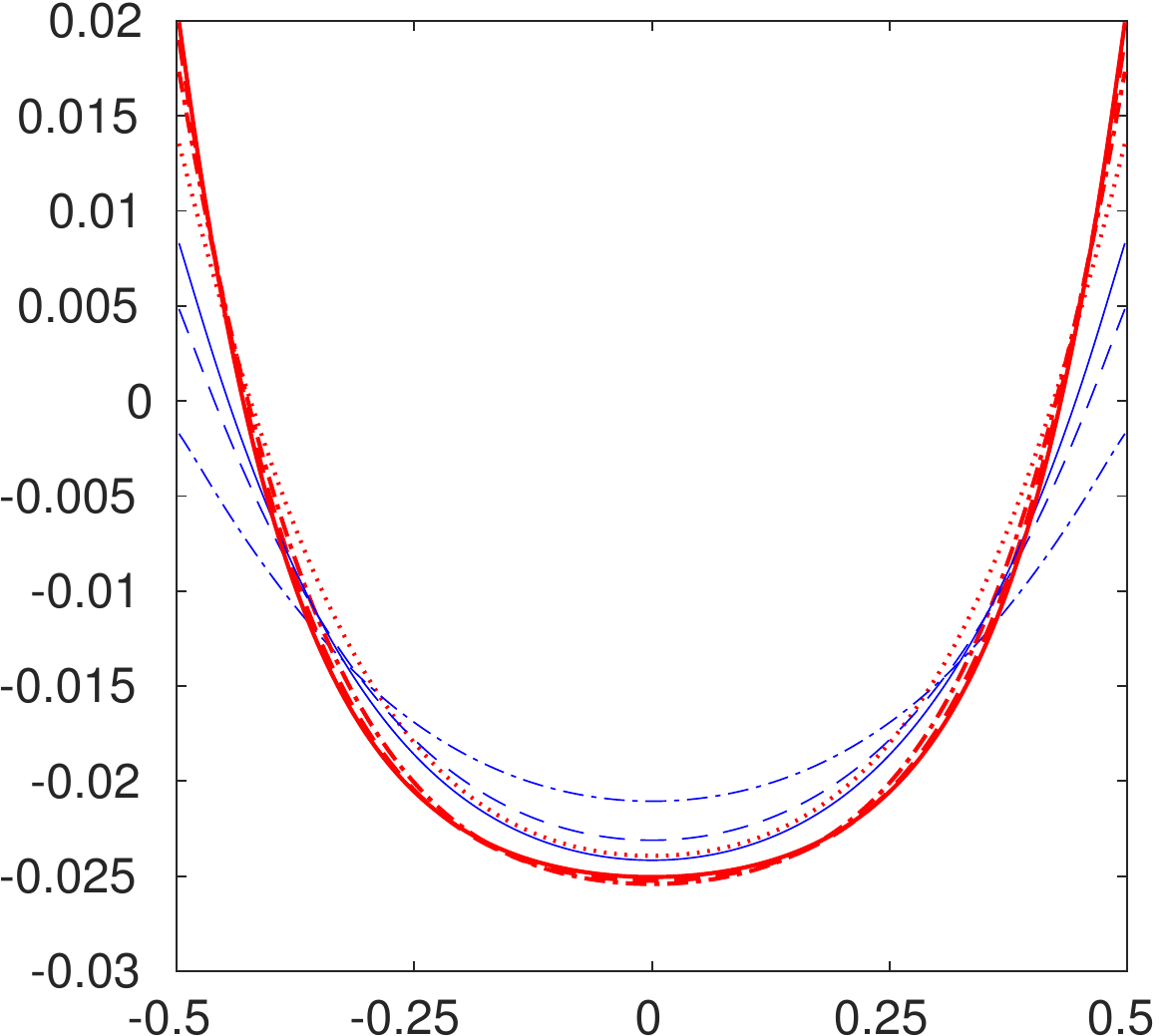}} 
\caption{Solution of the force driven Poiseuille flow on the uniform
  grid with $N=200$.}
  \label{fig:poiseuille-sol}
\end{figure}

With these settings, the steady-state solutions for density $\rho$,
temperature $\theta$, normal stress $\sigma_{11}$ and heat flux
$q_{2}$, obtained by the NMLM solver on the uniform grid with $N=200$,
are shown in \figurename~\ref{fig:poiseuille-sol}, which coincide well
with the steady-state solutions presented in \cite{hu2014nmg}, where
the first-order spatial discretization with $N=2048$ is employed.

\begin{figure}[!htb]
  \centering
  \subfloat[Poiseuille flow]{\label{fig:poiseuille-fourier:poiseuille}\includegraphics[width=0.49\textwidth]{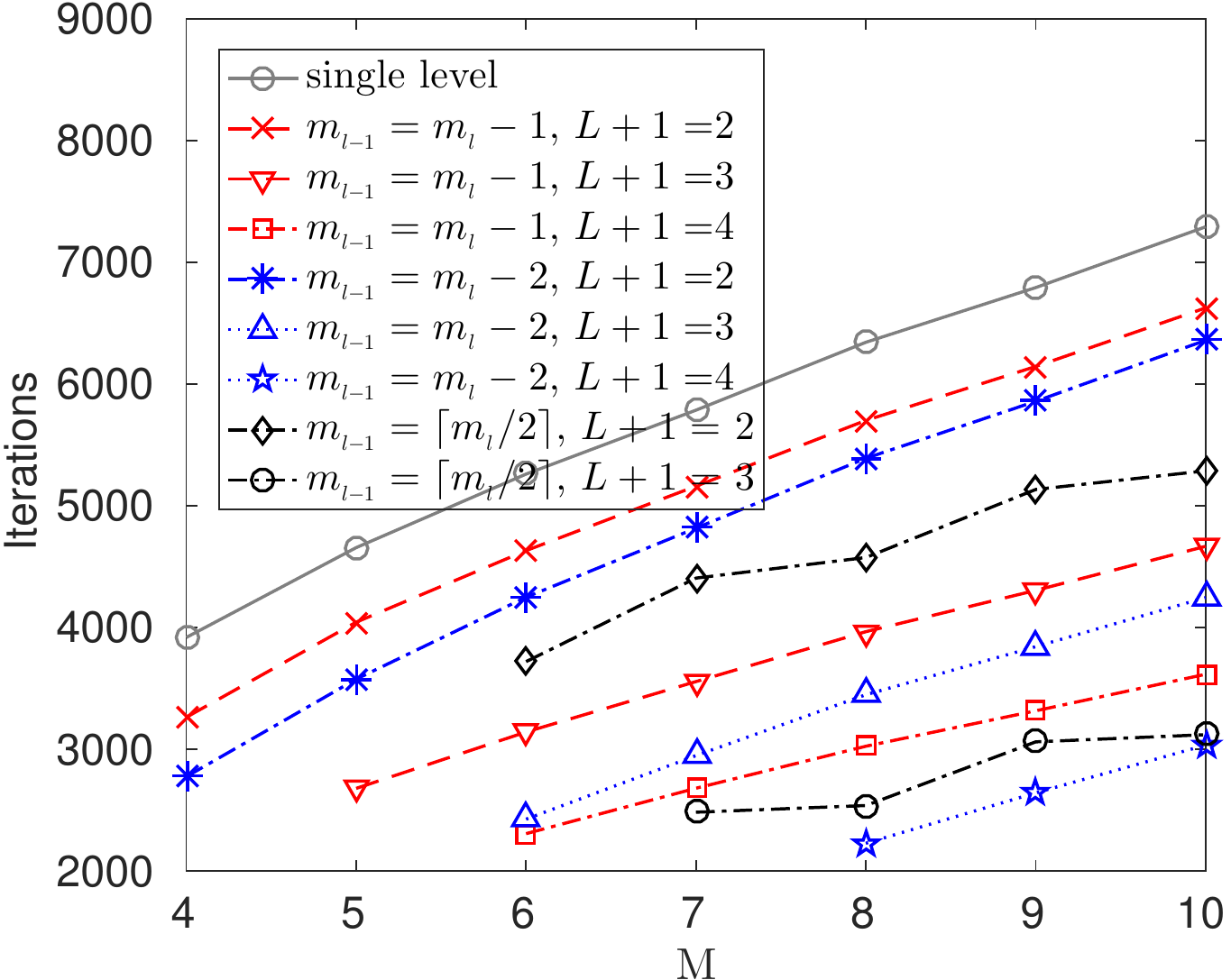}}\hfill
  \subfloat[Fourier flow]{\label{fig:poiseuille-fourier:fourier}\includegraphics[width=0.49\textwidth]{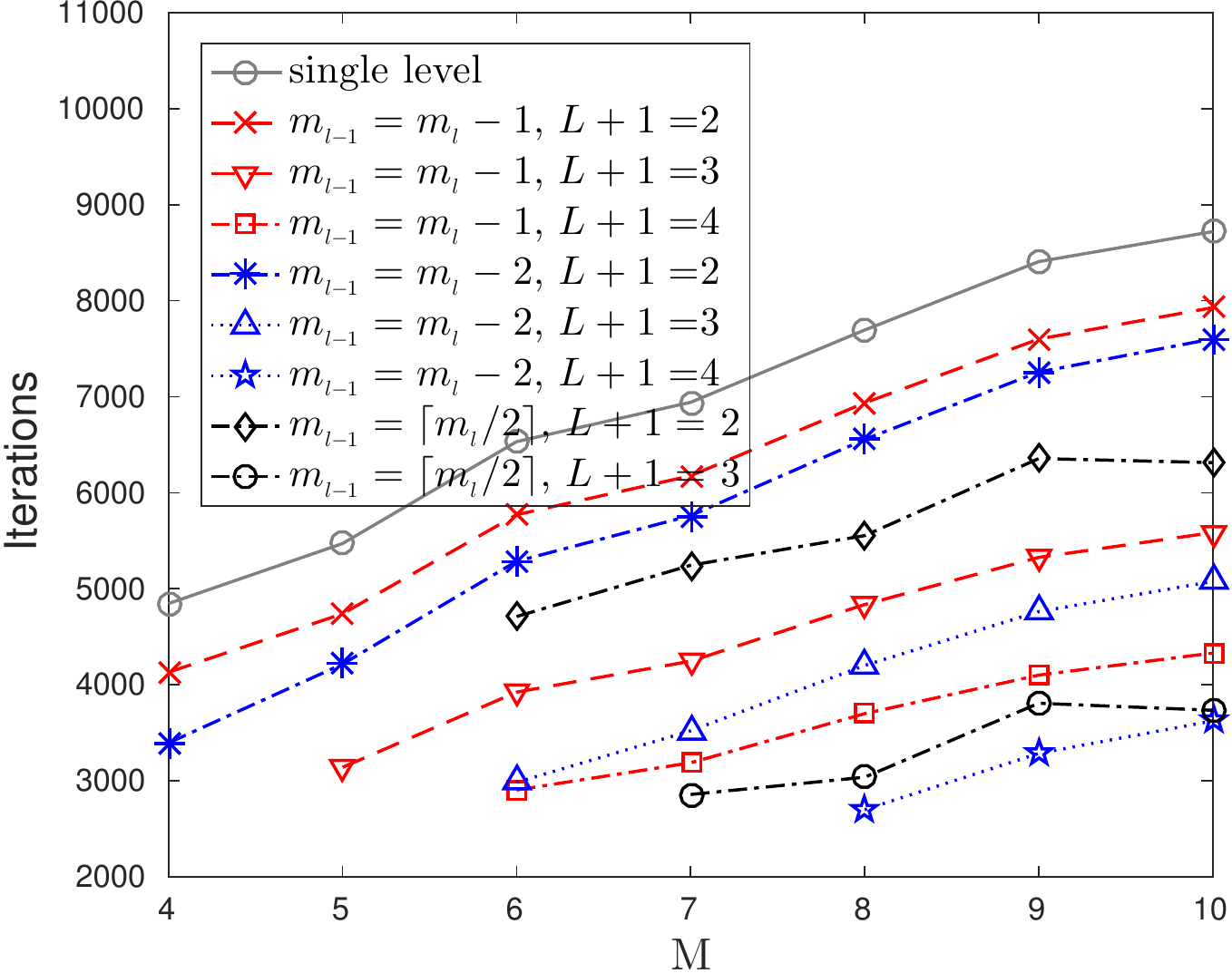}}
  \caption{Total number of iterations in terms of $M$ of the NMLM
    solver on the uniform grid of $N=200$. The total number of
    iterations of the single level solver is rescaled by a factor of
    $8$.}
  \label{fig:poiseuille-fourier:nmm-iters-orders}
\end{figure}

For the efficiency and behavior of the proposed NMLM solver, the tests
with various levels and order reduction strategies are performed on
three uniform grids with $N=100$, $200$ and $400$ for the moment model
with the order from $M=4$ to $10$.
As the Couette flow, only partial numerical results are presented
here. Specifically, the total number of iterations and the elapsed CPU
seconds are given in \tablename~\ref{tab:poiseuille:nmm-M45} for
$M=4,5$ and in
\tablename~\ref{tab:poiseuille:nmm-M10-p1}-\ref{tab:poiseuille:nmm-M10-p2}
for $M=10$ respectively. The corresponding convergence histories of
the tests on the uniform grid with $N=200$ are displayed in
\figurename~\ref{fig:poiseuille:nmm-res-history-N200-M4}-\ref{fig:poiseuille:nmm-res-history-N200-M10}.
The total number of iterations in terms of
$M$ is presented in
\figurename~\ref{fig:poiseuille-fourier:nmm-iters-orders}\subref{fig:poiseuille-fourier:poiseuille}.
All results show similar features as the tests of the Couette flow in
the case of $\Kn=0.1199$, which indicates the effectiveness of
the NMLM solver in accelerating the steady-state computation.

\begin{table}[!ht]
  \centering\footnotesize
  \begin{tabular}{c|c||c|c|c||c|cc|c}
    \hline\hline
    \multicolumn{2}{c||}{} & \multicolumn{3}{c||}{$M=4$} & \multicolumn{4}{c}{$M=5$}\\
    \hline
    \multicolumn{2}{c||}{} & & \multicolumn{1}{c|}{$m_{\sss l-1} = m_{\sss l} - 1$} & \multicolumn{1}{c||}{$m_{\sss l-1} = m_{\sss l} - 2$} & & \multicolumn{2}{c|}{$m_{\sss l-1} = m_{\sss l} - 1$} & \multicolumn{1}{c}{$m_{\sss l-1} = m_{\sss l} - 2$}\\ 
    \hline
    \multicolumn{2}{c||}{$L+1$} & 1 & 2 & 2 & 1 & 2 & 3 & 2 \\
    \hline
    \multirow{4}{*}{\begin{sideways}$N=100$\end{sideways}} 
    & $K$          &       15676&        1632&        1388&       18621&        2022&        1340&        1788\\
    & $T$          &     157.124&      87.112&      60.254&     268.329&     178.323&     140.127&     130.205\\
    & $K_s/K$      &       1.000&       9.605&      11.294&       1.000&       9.209&      13.896&      10.414\\
    & $T_s/T$      &       1.000&       1.804&       2.608&       1.000&       1.505&       1.915&       2.061\\
    \hline
    \multirow{4}{*}{\begin{sideways}$N=200$\end{sideways}} 
    & $K$          &       31349&        3263&        2775&       37251&        4044&        2679&        3576\\
    & $T$          &     566.855&     345.652&     239.170&    1068.763&     704.765&     559.148&     517.321\\
    & $K_s/K$      &       1.000&       9.607&      11.297&       1.000&       9.211&      13.905&      10.417\\
    & $T_s/T$      &       1.000&       1.640&       2.370&       1.000&       1.516&       1.911&       2.066\\
    \hline
    \multirow{4}{*}{\begin{sideways}$N=400$\end{sideways}} 
    & $K$          &       62694&        6524&        5590&       74509&        8088&        5358&        7152\\
    & $T$          &    2280.702&    1382.533&     971.918&    4291.040&    2818.351&    2232.552&    2070.655\\
    & $K_s/K$      &       1.000&       9.610&      11.215&       1.000&       9.212&      13.906&      10.418\\
    & $T_s/T$      &       1.000&       1.650&       2.347&       1.000&       1.523&       1.922&       2.072\\
    \hline\hline
  \end{tabular}
  \caption{Performance of the NMLM solver for the Poiseuille flow with $M=4,5$.}
  \label{tab:poiseuille:nmm-M45}
\end{table}

\begin{table}[!ht]
  \centering\footnotesize
  \begin{tabular}{c|c||ccccccc}
    \hline\hline
    \multicolumn{2}{c||}{} & \multicolumn{6}{c}{$m_{\sss l-1} = m_{\sss l} - 1$}\\
    \hline
    \multicolumn{2}{c||}{$L+1$} & 2 & 3 & 4 & 5 & 6 & 7 & 8 \\
    \hline
    \multirow{4}{*}{\begin{sideways}$N=100$\end{sideways}} 
    & $K$          &        3312&        2333&        1808&        1477&        1244&        1067&         921\\
    & $T$          &    1849.890&    1662.403&    1491.818&    1353.205&    1145.176&    1043.163&     905.261\\
    & $K_s/K$      &       8.803&      12.496&      16.125&      19.739&      23.436&      27.323&      31.655\\
    & $T_s/T$      &       1.388&       1.545&       1.721&       1.898&       2.243&       2.462&       2.837\\
    \hline
    \multirow{4}{*}{\begin{sideways}$N=200$\end{sideways}} 
    & $K$          &        6625&        4667&        3616&        2953&        2487&        2133&        1840\\
    & $T$          &    7498.892&    6641.554&    5983.719&    5250.433&    4644.889&    4111.988&    3637.646\\
    & $K_s/K$      &       8.804&      12.498&      16.131&      19.752&      23.454&      27.346&      31.701\\
    & $T_s/T$      &       1.304&       1.472&       1.634&       1.862&       2.105&       2.378&       2.688\\
    \hline
    \multirow{4}{*}{\begin{sideways}$N=400$\end{sideways}} 
    & $K$          &       13251&        9333&        7232&        5905&        4973&        4264&        3677\\
    & $T$          &   29539.018&   26645.046&   23642.063&   21103.438&   18830.564&   16646.179&   14537.716\\
    & $K_s/K$      &       8.804&      12.501&      16.132&      19.757&      23.460&      27.361&      31.729\\
    & $T_s/T$      &       1.301&       1.443&       1.626&       1.822&       2.042&       2.309&       2.644\\
    \hline\hline
  \end{tabular}
  \caption{Performance of the NMLM solver for the Poiseuille flow with $M=10$ (part I).}
  \label{tab:poiseuille:nmm-M10-p1}
\end{table}

\begin{table}[!ht]
  \centering\footnotesize
  \begin{tabular}{c|c||cccc|cc|c}
    \hline\hline
    \multicolumn{2}{c||}{} & \multicolumn{4}{c|}{$m_{\sss l-1} = m_{\sss l} - 2$} & \multicolumn{2}{c}{$m_{\sss l-1} = \lceil m_{\sss l} / 2 \rceil$} & \\
    \hline
    \multicolumn{2}{c||}{$L+1$} & 2 & 3 & 4 & 5 & 2 & 3 & 1\\
    \hline
    \multirow{4}{*}{\begin{sideways}$N=100$\end{sideways}} 
    & $K$          &        3179&        2123&        1517&        1059&        2643&        1560&       29154\\
    & $T$          &    1561.039&    1201.458&     888.561&     638.520&     989.148&     613.051&    2568.151\\
    & $K_s/K$      &       9.171&      13.732&      19.218&      27.530&      11.031&      18.688&       1.000\\
    & $T_s/T$      &       1.645&       2.138&       2.890&       4.022&       2.596&       4.189&       1.000\\
    \hline
    \multirow{4}{*}{\begin{sideways}$N=200$\end{sideways}} 
    & $K$          &        6359&        4246&        3034&        2116&        5287&        3120&       58329\\
    & $T$          &    6339.990&    4754.203&    3596.238&    2505.962&    3942.316&    2424.629&    9778.709\\
    & $K_s/K$      &       9.173&      13.737&      19.225&      27.566&      11.033&      18.695&       1.000\\
    & $T_s/T$      &       1.542&       2.057&       2.719&       3.902&       2.480&       4.033&       1.000\\
    \hline
    \multirow{4}{*}{\begin{sideways}$N=400$\end{sideways}} 
    & $K$          &       12718&        8491&        6067&        4231&       10574&        6239&      116668\\
    & $T$          &   25255.174&   18999.092&   14340.679&   10061.378&   15946.033&    9624.905&   38442.738\\
    & $K_s/K$      &       9.173&      13.740&      19.230&      27.575&      11.033&      18.700&       1.000\\
    & $T_s/T$      &       1.522&       2.023&       2.681&       3.821&       2.411&       3.994&       1.000\\
    \hline\hline
  \end{tabular}
  \caption{Performance of the NMLM solver for the Poiseuille flow with $M=10$ (part II).}
  \label{tab:poiseuille:nmm-M10-p2}
\end{table}

\newcommand\drawPoiseuilleNMMHistory[2]{\begin{figure}[!htb]
  \centering
  {\includegraphics[width=0.5\textwidth]{poiseuille_nmm_res_iters_N#1_M#2.pdf}}\hfill
  {\includegraphics[width=0.5\textwidth]{poiseuille_nmm_res_cputime_N#1_M#2.pdf}}
  \caption{Convergence history of the NMLM solver for the Poiseuille
    flow with $M=#2$ on the uniform grid of $N=#1$.}
  \label{fig:poiseuille:nmm-res-history-N#1-M#2}
\end{figure}
}
\drawPoiseuilleNMMHistory{200}{4}
\drawPoiseuilleNMMHistory{200}{5}
\drawPoiseuilleNMMHistory{200}{10}

\subsection{The Fourier flow}
\label{sec:num-ex-fourier}
\begin{figure}[!htb]
  \centering
  \subfloat[Density, $\rho$]{\includegraphics[width=0.49\textwidth]{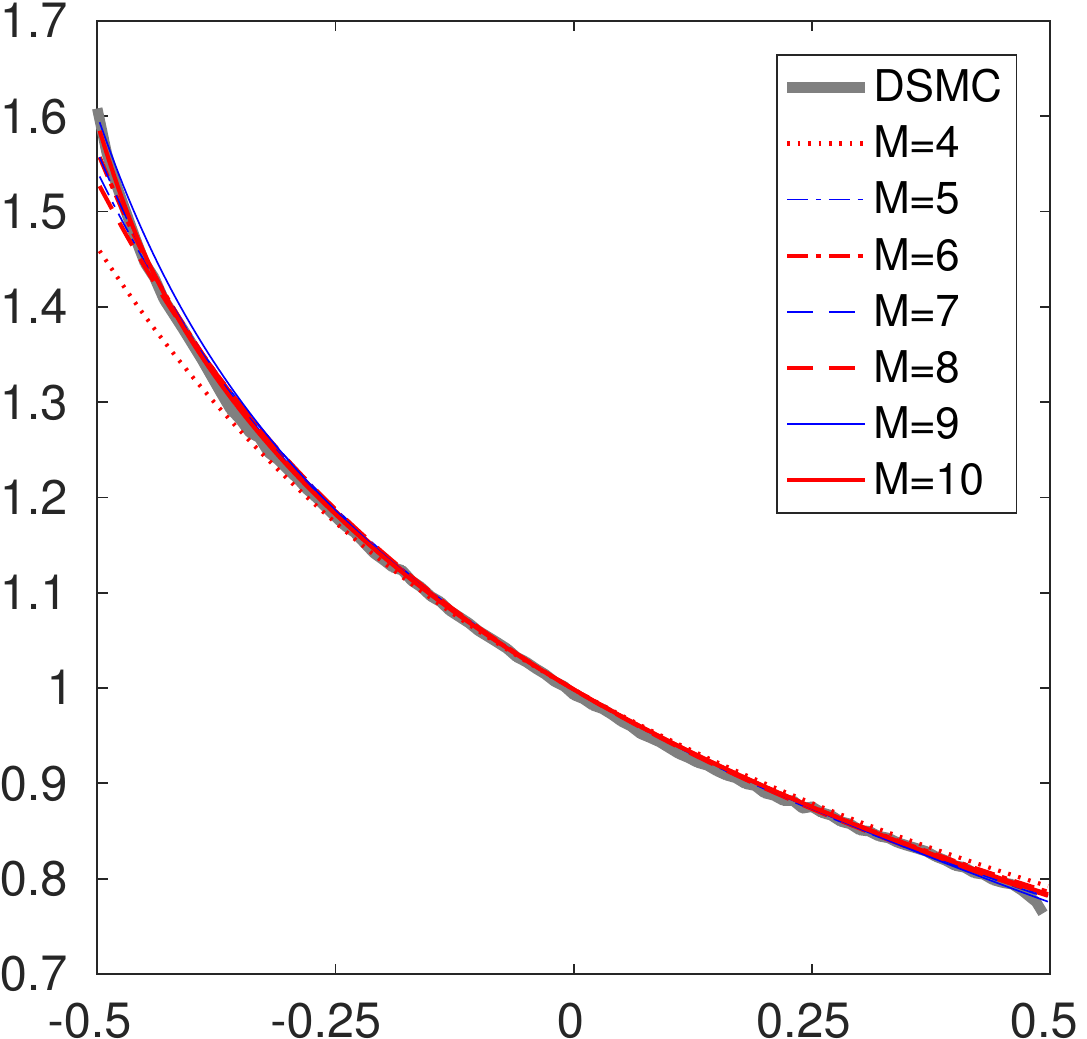}}\hfill
  \subfloat[Temperature, $\theta$]{\includegraphics[width=0.49\textwidth]{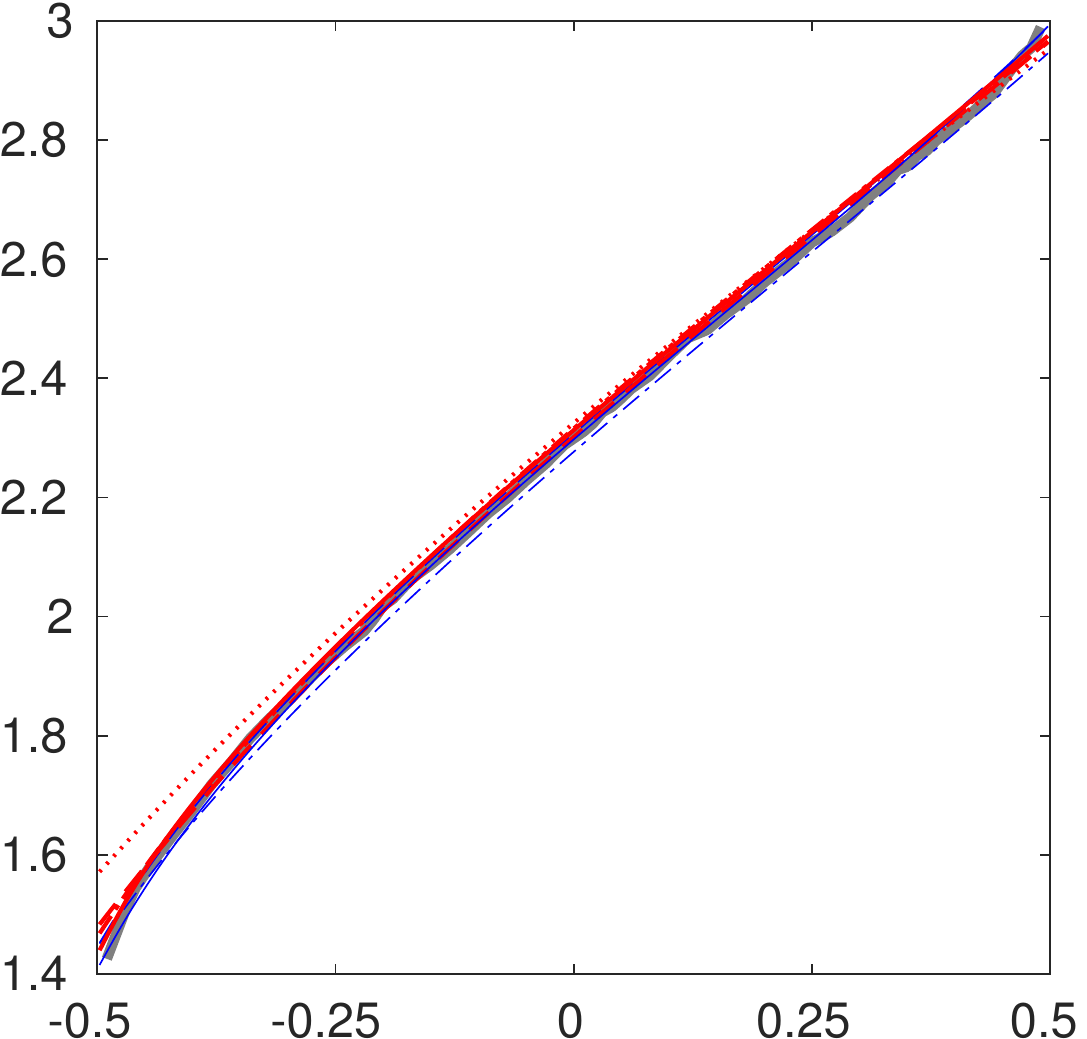}}
  \caption{Solution of the Fourier flow on the uniform grid with
    $N=200$.}
  \label{fig:fourier-sol}
\end{figure}

The last benchmark test is the Fourier flow which also investigates
the motion of the gas between two infinite parallel plates with a
distance of $L_{D}=1$. In contrast to the previous examples, both
plates are stationary, while their temperatures are different. The gas
is driven by the difference of temperatures between the two plates,
and could reach a steady state in the absence of external force, that
is, $\bF\equiv 0$. To reproduce the results in \cite{Microflows1D,
  Wadsworth}, the gas of helium with the viscosity index $w = 0.657$
and the Knudsen number $\Kn=0.1044$ for the collision frequency
\eqref{eq:vhs-nu} is considered. The temperature on the left plate and
the right plate are set to be 0.2894 and 1.0769
respectively. Numerical solutions for density $\rho$ and temperature
$\theta$, obtained by the NMLM solver on the uniform grid with
$N=200$, are shown in \figurename~\ref{fig:fourier-sol}. The solutions
obtained by the DSMC (Direct Simulation of Monte Carlo) method
\cite{Wadsworth} are provided as a reference. It can be observed that
the solutions of the moment model converge and match the DSMC solution
well as the order $M$ increases.

As for the performance of the NMLM solver, the tests with various
levels and order reduction strategies are also performed on three
uniform grids with $N=100$, $200$ and $400$ for the moment model with
the order from $M=4$ to $10$. Due to the same reason, only partial
numerical results are presented here. That is, the total number of
iterations and the elapsed CPU seconds for $M=10$ are given in
\tablename~\ref{tab:fourier:nmm-M10-p1}-\ref{tab:fourier:nmm-M10-p2}. The
corresponding convergence histories of the tests on the uniform grid
with $N=200$ are displayed in
\figurename~\ref{fig:fourier:nmm-res-history-N200-M10}. And the total
number of iterations in terms of $M$ is plotted in
\figurename~\ref{fig:poiseuille-fourier:nmm-iters-orders}\subref{fig:poiseuille-fourier:fourier}.
Again, all results show similar features as the tests of the Couette
flow in the case of $\Kn=0.1199$ and the tests of the Poiseuille flow.
Therefore, the proposed NMLM solver is indeed able to accelerate the
steady-state computation significantly.

\begin{table}[!ht]
  \centering\footnotesize
  \begin{tabular}{c|c||ccccccc}
    \hline\hline
    \multicolumn{2}{c||}{} & \multicolumn{6}{c}{$m_{\sss l-1} = m_{\sss l} - 1$}\\
    \hline
    \multicolumn{2}{c||}{$L+1$} & 2 & 3 & 4 & 5 & 6 & 7 & 8 \\
    \hline\hline
    \multirow{4}{*}{\begin{sideways}$N=100$\end{sideways}} 
    & $K$          &        3965&        2794&        2167&        1771&        1493&        1282&        1109\\
    & $T$          &    1860.858&    1694.069&    1523.898&    1363.943&    1209.693&    1075.778&     930.812\\
    & $K_s/K$      &       8.799&      12.486&      16.099&      19.699&      23.367&      27.213&      31.458\\
    & $T_s/T$      &       1.592&       1.749&       1.944&       2.172&       2.449&       2.754&       3.183\\
    \hline
    \multirow{4}{*}{\begin{sideways}$N=200$\end{sideways}} 
    & $K$          &        7928&        5585&        4331&        3539&        2983&        2560&        2214\\
    & $T$          &    7415.489&    6720.777&    6035.377&    5384.300&    4840.328&    4236.370&    3750.768\\
    & $K_s/K$      &       8.799&      12.490&      16.106&      19.711&      23.385&      27.248&      31.507\\
    & $T_s/T$      &       1.593&       1.758&       1.958&       2.194&       2.441&       2.789&       3.150\\
    \hline
    \multirow{4}{*}{\begin{sideways}$N=400$\end{sideways}} 
    & $K$          &       15851&       11167&        8659&        7074&        5962&        5116&        4424\\
    & $T$          &   29739.900&   26760.605&   24140.459&   21534.845&   19151.365&   17020.835&   14997.699\\
    & $K_s/K$      &       8.800&      12.491&      16.109&      19.718&      23.396&      27.265&      31.530\\
    & $T_s/T$      &       1.527&       1.697&       1.881&       2.108&       2.371&       2.667&       3.027\\
    \hline\hline
  \end{tabular}
  \caption{Performance of the NMLM solver for the Fourier flow with $M=10$ (part I).}
  \label{tab:fourier:nmm-M10-p1}
\end{table}

\begin{table}[!ht]
  \centering\footnotesize
  \begin{tabular}{c|c||cccc|cc|c}
    \hline\hline
    \multicolumn{2}{c||}{} & \multicolumn{4}{c|}{$m_{\sss l-1} = m_{\sss l} - 2$} & \multicolumn{2}{c}{$m_{\sss l-1} = \lceil m_{\sss l} / 2 \rceil$} & \\
    \hline
    \multicolumn{2}{c||}{$L+1$} & 2 & 3 & 4 & 5 & 2 & 3 & 1\\
    \hline
    \multirow{4}{*}{\begin{sideways}$N=100$\end{sideways}} 
    & $K$          &        3803&        2540&        1814&        1252&        3159&        1870&       34887\\
    & $T$          &    1579.625&    1213.639&     907.115&     631.572&     975.638&     603.752&    2962.720\\
    & $K_s/K$      &       9.174&      13.735&      19.232&      27.865&      11.044&      18.656&       1.000\\
    & $T_s/T$      &       1.876&       2.441&       3.266&       4.691&       3.037&       4.907&       1.000\\
    \hline
    \multirow{4}{*}{\begin{sideways}$N=200$\end{sideways}} 
    & $K$          &        7603&        5076&        3624&        2502&        6314&        3735&       69756\\
    & $T$          &    6361.864&    4843.804&    3627.384&    2553.007&    3926.485&    2419.668&   11815.652\\
    & $K_s/K$      &       9.175&      13.742&      19.248&      27.880&      11.048&      18.676&       1.000\\
    & $T_s/T$      &       1.857&       2.439&       3.257&       4.628&       3.009&       4.883&       1.000\\
    \hline
    \multirow{4}{*}{\begin{sideways}$N=400$\end{sideways}} 
    & $K$          &       15203&       10149&        7244&        5001&       12622&        7465&      139487\\
    & $T$          &   25136.041&   19347.345&   14581.706&   10236.093&   15873.863&    9702.212&   45400.647\\
    & $K_s/K$      &       9.175&      13.744&      19.256&      27.892&      11.051&      18.685&       1.000\\
    & $T_s/T$      &       1.806&       2.347&       3.114&       4.435&       2.860&       4.679&       1.000\\
    \hline\hline
  \end{tabular}
  \caption{Performance of the NMLM solver for the Fourier flow with $M=10$ (part II).}
  \label{tab:fourier:nmm-M10-p2}
\end{table}

\newcommand\drawFourierNMMHistory[2]{\begin{figure}[!htb]
  \centering
  {\includegraphics[width=0.5\textwidth]{fourier_nmm_res_iters_N#1_M#2.pdf}}\hfill
  {\includegraphics[width=0.5\textwidth]{fourier_nmm_res_cputime_N#1_M#2.pdf}}
  \caption{Convergence history of the NMLM solver for the Fourier
    flow with $M=#2$ on the uniform grid of $N=#1$.}
  \label{fig:fourier:nmm-res-history-N#1-M#2}
\end{figure}
}
\drawFourierNMMHistory{200}{10}


\section{Concluding remarks}
\label{sec:conclusion}

A steady-state solver for microflows with high-order
  moment model was successfully proposed in this paper, which significantly
  improved the efficiency of the one in \cite{hu2016acceleration} from the
  following approaches:
\begin{itemize}
\item Linear reconstruction is adopted for high-resolution spatial
  discretization, so that remarkable reduction for degrees of freedom
  in spatial space is obtained without loss of accuracy.
\item A relaxation parameter is introduced in the correction step to
  enhance the stability of the solver such that more levels can be
  applied in the solver. 
\item The computation of the correction step is also simplified a lot
  in comparison to the way used in \cite{hu2016acceleration}.
\item Heun's method is taken as the smoother in each level to further
  improve the robustness of the NMLM solver in the situation when many
  levels are involved.
\end{itemize}
The performance of the new NMLM solver is numerically investigated by
three benchmark problems in microflows. Various order reduction
strategies for the choice of the order sequence of the NMLM solver
have been tested. For each order reduction strategy, the convergence
rate of the resulting NMLM solver is improved as the total levels
increases. Among these order reduction strategies, it is shown that
the most efficient strategy is
$m_{\sss l-1} = \lceil m_{\sss l} / 2 \rceil$, and the second strategy
is $m_{\sss l-1} = m_{\sss l} - 2$. In summary, it is demonstrated
that the new NMLM solver can further improve the efficiency of
steady-state computations even for the moment model with a relatively
small order, such as $M=4$ and $5$. As a result, the idea of using the
lower-order moment model correction is very promising to accelerate
the steady-state simulation, and may also be valuable for problems
described by other hierarchical models.

Additionally, the NMLM solver behaves similarly to the single level
solver, as the order $M$ or the spatial grid number $N$
increases. Research works on combination of the lower-order moment
model correction with the spatial coarse grid correction are ongoing.

\section*{Acknowledgements} The authors would like to thank Prof. Ruo Li at
Peking University, China for the constructive suggestions to this work. The
research of Zhicheng Hu is partially supported by the National Natural Science
Foundation of China (11601229), and the Natural Science Foundation of Jiangsu
Province of China (BK20160784). The research of Guanghui Hu is partially
supported by the FDCT of Macao SAR (029/2016/A1), the MYRG of University of
Macau (MYRG2017-00189-FST), and the National Natural Science Foundation of China
(11401608).



\bibliographystyle{plain}
\bibliography{ref/steadystate,ref/article,ref/tiao}
\end{document}